\documentclass[journal]{IEEEtran}
\ifCLASSINFOpdf
\else
\fi
\usepackage{amsmath}
\usepackage{amsthm}
\usepackage{color}
\usepackage{pgfplots} 
\usepackage{pgf}
\usepackage{tikz}
\usetikzlibrary{arrows,automata,patterns,chains}
\usepackage{amssymb}
\usepackage{mathtools}
\usepackage{algorithm}
\usepackage{enumitem}
\usepackage{textcomp}
\usepackage{algpseudocode}
\usepackage{amsfonts}
\usepackage{subfiles}
\usepackage{float}
\usepackage{tkz-graph}

\makeatletter
\renewcommand*\env@matrix[1][*\c@MaxMatrixCols c]{%
  \hskip -\arraycolsep
  \let\@ifnextchar\new@ifnextchar
  \array{#1}}
\makeatother  

\newtheorem{assumption}{Assumption}
\newtheorem{theorem}{Theorem}
\newtheorem{lemma}{Lemma}

\newfloat{algorithm}{t}{lop}
\newcommand{\subscript}[2]{$#1 _ #2$}

\newcommand{\pgfplotsdrawaxis}{\pgfplots@draw@axis}
\makeatother
\pgfplotsset{only axis on top/.style={axis on top=false, after end axis/.code={
             \pgfplotsset{axis line style=opaque, ticklabel style=opaque, tick style=opaque,
                          grid=none}\pgfplotsdrawaxis}}}
\newcommand{\drawge}{-- (rel axis cs:1,0) -- (rel axis cs:1,1) -- (rel axis cs:0,1) \closedcycle}
\newcommand{\drawle}{-- (rel axis cs:1,1) -- (rel axis cs:1,0) -- (rel axis cs:0,0) \closedcycle}
\begin{document}
\title{A geometrically converging dual method for distributed optimization over time-varying graphs}

\author{Marie~Maros,~\IEEEmembership{Student Member,~IEEE,}
        and Joakim~Jald\'{e}n,~\IEEEmembership{Member,~IEEE,}% <-this % stops a space

\thanks{This project has received funding from the European Research Council (ERC) under the European Union's Horizon 2020 research and innovation programme (grant agreement No 742648)}
\thanks{This paper will be presented in part at the 57th IEEE Conference on Decision and Control, Florida, USA, December 2018.}
\thanks{The authors are with the department of Information Science and Engineering, School of Electrical Engineering and Computer Science, KTH, Sweden. (e-mails:\emph{mmaros@kth.se}, \emph{jalden@kth.se})}% <-this % stops a space
}


\markboth{TAC}%
{Shell \MakeLowercase{\textit{et al.}}: Bare Demo of IEEEtran.cls for IEEE Journals}



\maketitle

\begin{abstract}
In this paper we consider a distributed convex
optimization problem over time-varying undirected networks. We propose a dual method, primarily averaged network dual ascent (PANDA), that is proven to converge R-linearly to the optimal point
given that the agents’ objective functions are strongly convex
and have Lipschitz continuous gradients. Like dual decomposition, PANDA
requires half the amount of variable exchanges per iterate
of methods based on DIGing, and can provide with practical improved
performance as empirically demonstrated.
\end{abstract}


% \fi
\IEEEpeerreviewmaketitle

\section{Introduction}
\IEEEPARstart{M}{any} engineering applications rely on effectively solving problems in a distributed manner. These include applications in distributed estimation and control \cite{xiao,application,formation,estimation}, decentralized source localization \cite{localization}, power grids \cite{giannakis,grid} and distributed learning \cite{learning, learning2, social_learning}. In these set-ups agents have partial information of the problem they are to cooperatively solve, and gathering all the information in a fusion center may not be affordable nor desirable. The agents may have gathered massive amounts of data making its exchange expensive, or the agents may be connected via a wireless network resulting in communication costs. Due to the random nature of the wireless channel the communication links can be unreliable. This justifies the need for procedures according to which the agents cooperate to solve the overall optimization problem via an unreliable network. 
In this paper w therefore focus on solving  the problem
\begin{equation}
\label{eq:original}
\underset{\bar{\mathbf{x}} \in \mathbb{R}^p}{\text{min}} \quad \bar{f}(\bar{\mathbf{x}}) \triangleq \sum_{i=1}^n f_i(\bar{\mathbf{x}}),
\end{equation}
in a decentralized manner by having $n$ agents cooperate over a time-varying network. We will represent the network with a sequence of graphs. Agent $i$ is only aware of its individual objective function $f_i: \mathbb{R}^p \to \mathbb{R}$ and the goal is that the agents  cooperatively solve \eqref{eq:original} by exchanging information without explicitly exchanging their individual objective functions. 
% Distributed optimization: undirected, then directed

Many algorithms to solve \eqref{eq:original} in a distributed manner over a static and undirected network have been proposed. First order primal methods such as decentralized (sub)-gradient descent typically require a diminishing step-size to converge \cite{decentralized_gradient,subgradient,nedic}. A diminishing step-size will imply slower convergence rates. If the individual objective functions are strongly convex and have Lipschitz continuous gradients, decentralized gradient descent will converge linearly to a neighborhood of the optimal solution when using a fixed step-size. However, the size of the neighborhood will be proportional to the step-size leading to a speed-accuracy dilemma. In order to circumvent this problem, algorithms that incorporate an error correcting scheme to the iterates of decentralized gradient descent, such as the decentralized exact first-order algorithm (EXTRA) \cite{extra} and \cite{harnesing} have been proposed. The incorporation of the error correcting schemes allows the algorithms to converge linearly to the optimal point with constant, but appropriately chosen, step-sizes. Proximal versions of such algorithms have also been proposed \cite{pextra}. The Directed EXTRA (DEXTRA) \cite{dextra} algorithm was proposed as an extension of EXTRA \cite{extra} to directed graphs. While DEXTRA preserves the convergence properties of EXTRA it requires the step-size to fall within a strictly positive interval. Further, this interval is not always non-empty \cite{dextra}. The Accelerated Distributed Direct Optimization (ADD-OPT), an alternative to DEXTRA without the mentioned issues, was recently proposed in \cite{addopt}. Further extensions include \cite{frost} in which the step-sizes of different agents need not be identical. However, extensions of \cite{extra,harnesing} and \cite{pextra} to time-varying networks do not exist in the literature.

While the methods in \cite{extra,harnesing,pextra} were proposed as modified versions of gradient descent they can also be interpreted as primal-dual methods as explained in \cite{interpret}. Distributed Dual and Primal-Dual methods typically converge linearly to the optimal point using a constant step-size \cite{admm, dadmm, dlm, extra, harnesing, pextra, dual_decomposition, tutorial, optimal}. The distributed alternating direction method of multipliers (ADMM) \cite{admm}  has been shown to converge linearly over an undirected and static graph if the individual objective functions are strongly convex and have Lipschitz continuous gradients. Further, ADMM is also shown to converge Q-linearly when the sum of the obje tive functions is strongly convex but the individual objective functions are not \cite{maros}. \cite{random_networks} proposed an ADMM-based method that converges sub-linearly over random time-varying undirected graphs. Further, \cite{fenchel} established that Dual Decomposition \cite{dual_decomposition} converges to the primal optimal solution at rate $\mathcal{O}(\frac{1}{\sqrt{k}})$ on time-varying undirected graphs given that the individual objective functions are strongly convex. By proposing a different error correction scheme than \cite{extra,harnesing} the authors of \cite{interpret} proposed decentralized inexact gradient tracking (DIGing), a method that converges linearly to the optimal solution even if the undirected graph is time-varying. Further, in \cite{interpret} the authors also develop a version of DIGing, Push-DIGing, capable of converging linearly over directed time-varying graphs. This is achieved by combining DIGing and the push-sum protocol. In-network non-convex optimization (NEXT) \cite{next} is a method to solve non-convex problems of the form \eqref{eq:original} on undirected graphs. The authors in \cite{non_convex} proposed the successive convex approximation over time-varying directed graphs (SONATA) to solve non-convex optimization problems over time-varying directed graphs. SONATA can be particularized to yield the same iterates as Push-DIGing. Under the assumption that the individual functions $f_i$ are strongly convex and have Lipschitz continuous gradients NEXT, SONATA, DIGing and Push-DIGing will converge linearly to the optimal point. Further, the four of them are primal methods that rely on tracking the average primal variable and the average gradient. This procedure, while effective, requires that the agents exchange the primal variables and the gradients at each iterate.

In this paper, we propose \underline{p}rimarily \underline{a}veraged \underline{n}etwork \underline{d}ual \underline{a}scent (PANDA). Different to dual decomposition \cite{dual_decomposition} PANDA actively tracks the average of the primal variable so as to use it as the dual gradient.  Assuming that the objective functions are strongly convex and have Lipschitz continuous gradients, we establish that PANDA converges R-linearly over time-varying undirected graphs. To our knowledge this is the first dual method that can be shown to converge linearly over time-varying networks. Further, PANDA does not require the usage of symmetric mixing matrices as opposed to Dual Decomposition. This makes the possibility of tweaking the iterates in PANDA to deal with directed graphs more realistic. While the iterates of PANDA are computationally more expensive than those of DIGing PANDA does not require the exchange of gradients. In terms of convergence rate, PANDA performs better in practice as shown in our experiments.

The rest of the paper is structured as follows. In Section \ref{section:algorithm} we intuitively and formally introduce PANDA as an alternative to Dual Decomposition that can handle time-varying graphs. Further, we provide the formal statement regarding PANDA's convergence in the form of a convergence theorem. Section \ref{section:proof} is devoted to proving the theorem introduced in Section \ref{section:algorithm}. In particular, Section \ref{section:proof} is divided in two essentially different subsections. In the first section we establish a number of Lemmas that are required to be able to invoke the small gain theorem \cite{interpret}. In the following sub-section we use the small-gain theorem to finally PANDA's R-linear convergence, with inspiration from \cite{interpret}. Section \ref{section:numerical} is devoted to numerical experiments that illustrate the performance of PANDA as compared to other state of the art methods. Finally we end the paper with some concluding remarks.

\section{Algorithm \label{section:algorithm}}
In this section we intuitively introduce PANDA. Afterwards, we introduce the assumptions under which we formally provide PANDA's convergence guarantees.
\subsection{Intuition}
Let $\mathcal{G} = (\mathcal{V},\mathcal{E})$ denote the undirected graphs representing the network that connects the agents, indexed by $i=1,\hdots,n.$ The agents are represented by the graph's vertices $\mathcal{V}$ and the communication links by the edges $\mathcal{E}.$ If two agents $i$ and $j$ can communicate, i.e. if $(i,j)\in \mathcal{E},$ then $(j,i) \in \mathcal{E}$ since the graph is undirected. The set of agents with which $i$ can directly communicate is the neighborhood of $i$ and is denoted by $\mathcal{N}_i.$
In order to solve \eqref{eq:original} in a distributed manner over $\mathcal{G}$ it is generally required that an individual copy of the variable $\bar{\mathbf{x}}$ is assigned to each agent. Thus, we formulate the optimization problem
\begin{equation}
\label{eq:on_graph}
\underset{\mathbf{x} \in \mathbb{R}^{np}}{\min} \quad f(\mathbf{x}) \triangleq \sum_{i=1}^n f_i(\mathbf{x}_i), \quad \text{s.t. } \mathbf{x}_i = \mathbf{x}_j,\forall i\in \mathcal{V}, \forall j \in \mathcal{N}_i,
\end{equation}
where $\mathbf{x} \triangleq (\mathbf{x}_1^T,\hdots,\mathbf{x}_n^T)^T.$ The constraint  in \eqref{eq:on_graph} enforces $\mathbf{x}$ to be consensual, i.e. it enforces that all nodes agree on $\mathbf{x}_1 =\mathbf{x}_2 = \hdots= \mathbf{x}_n.$ Problems \eqref{eq:on_graph} and \eqref{eq:original} are equivalent in the sense that if $\bar{\mathbf{x}}^{\star}$ solves \eqref{eq:original}, then $\mathbf{x}^{\star} = (\bar{\mathbf{x}}^{\star},\hdots,\bar{\mathbf{x}}^{\star})$ solves \eqref{eq:on_graph} provided that the graph $\mathcal{G}$ is connected. In some methods, the constraint in \eqref{eq:on_graph} is treated implicitly via the averaging of the variables and gradients of the different agents \cite{extra,interpret}. Under the assumption that $\mathcal{G}$ is connected, dual or primal-dual methods formulate a constrained optimization problem such as the one in \eqref{eq:on_graph} which they will explicitly solve \cite{dual_decomposition,optimal,dadmm}. Note that in the formulation in \eqref{eq:on_graph} all agents implicitly treat all other agents equally. More generally, an agent $i$ may assign different importances to information coming from different agents within its neighborhood. Such an optimization problem can be expressed as
\begin{equation}
\label{eq:matrix}
\underset{\mathbf{x} \in \mathbb{R}^{np}}{\min} \quad f(\mathbf{x}), \quad \text{s.t. } \left(\mathbf{U}^{1/2}\otimes \mathbf{I}_p\right)\mathbf{x} = \mathbf{0},
\end{equation} 
where the matrix $\mathbf{U} \in \mathbb{R}^{n\times n}$ has the following properties \cite{optimal}:
\begin{enumerate}[label = (\subscript{P}{{\arabic*}})]
\item $\mathbf{U}$ is an $n \times n $ symmetric matrix,
\item $\mathbf{U}$ is positive semi-definite,
\item \label{property:consensual} $\text{null}\{\mathbf{U}\} = \text{span}\{\mathbf{1}_n\},$ where $\mathbf{1}_n = (1,\hdots,1)^T \in \mathbb{R}^n,$
\item $\mathbf{U}$ is defined on the edges of the graph, i.e. $u_{ij} \neq 0$ only if $(i,j) \in \mathcal{E}.$
\end{enumerate}
Note that property \ref{property:consensual} enforces that $\mathbf{x}$ is consensual. However, as noted, this property can only be fulfilled when the underlying graph $\mathcal{G}$ is connected. Problem \eqref{eq:matrix} can be solved in a distributed manner using dual ascent as follows
\begin{subequations}
\begin{align}
\mathbf{x}_i(k+1) &:= \text{arg } \underset{\mathbf{x}_i \in \mathbb{R}^{p}}{\min} \quad f_i(\mathbf{x}_i) - \mathbf{y}_i^T\mathbf{x}_i \\
\mathbf{y}_{i}(k+1) &:= \mathbf{y}_i(k) -  c \!\!\!\!\!\! \sum_{j \in \mathcal{N}_i \cup \{i\}}\!\!\!\!\!\! u_{ij}\mathbf{x}_j(k+1) \label{eq:dual_ascent_y}
\end{align}
\end{subequations}
where $c > 0$ is an appropriately selected step-size and $u_{ij}$ is the weight node $i$ assigns to the information coming from node $j.$ Note that the quantity $u_{ij}$ corresponds to element $(i,j)$ of the matrix $\mathbf{U}$. To illustrate the short-comings of dual ascent and to motivate the proposal of PANDA we will now analyze what happens when the network and consequently the graph change over time. To this end, let $\{\mathcal{G}(k)\}$ denote a sequence of graphs such that $\mathcal{G}(k) = (\mathcal{V},\mathcal{E}(k)).$ Consequently the iterate \eqref{eq:dual_ascent_y} becomes
\begin{equation}
\mathbf{y}_{i}(k+1) := \mathbf{y}_i(k) - c \!\!\!\!\!\!\!\!\sum_{j \in \mathcal{N}_i(k) \cup \{i\}}\!\!\!\!\!\!\!\! u_{ij}(k)\mathbf{x}_j(k+1),
\end{equation}
where both the neighborhood of $i,$ $\mathcal{N}_i(k),$ and the weights $u_{ij}(k)$ are now changing with time. Assume further, that at time $k$ the graph $\mathcal{G}(k)$ is not connected. This will imply several things. First of all, property \ref{property:consensual} is broken and therefore problems \eqref{eq:matrix} and \eqref{eq:original} are no longer equivalent. This is because the problem formulation in \eqref{eq:matrix} is graph dependent. Further, at time $k$ we are performing a step towards solving independent optimization problems within the connected sub-graphs instead of attempting to solve \eqref{eq:original}. The mentioned problems are the consequence of two essential issues: 
\begin{itemize}
\item the problem formulation \eqref{eq:matrix} is graph dependent,
\item there is no explicit averaging across time to take advantage of the different graph structures given by the sequence $\{\mathcal{G}(k)\}.$
\end{itemize}
Addressing the first issue can be done by identifying something that all connected graphs have in common. For this we write the dual function of \eqref{eq:matrix} as
\begin{equation}
\label{eq:dual}
\underset{\boldsymbol{\lambda} \in \mathbb{R}^{np}}{\min} \quad f^*\left(\left(\mathbf{U}^{1/2} \otimes \mathbf{I}_p\right)\boldsymbol{\lambda}\right),
\end{equation}
where $f^*$ denotes the convex conjugate
\begin{equation}
f^*(\mathbf{y}) = \underset{\mathbf{x}\in \mathbb{R}^{np}}{\sup} \quad f(\mathbf{x}) - \mathbf{y}^T\mathbf{x}.
\end{equation}
Problem \eqref{eq:dual} can be reformulated equivalently by adding a variable $\mathbf{y}$ as follows
\begin{equation}
\underset{\mathbf{y} \in \mathbb{R}^{np},\,\boldsymbol{\lambda} \in \mathbb{R}^{np}}{\min} \quad f^*(\mathbf{y}) \quad \text{s.t. } \left( \mathbf{U}^{1/2}\otimes \mathbf{I}_p \right) \boldsymbol{\lambda} = \mathbf{y}.
\end{equation}
Note that $\text{span}\{\mathbf{U}^{1/2}\} = \mathbb{R}^n \setminus \text{span}\{\mathbf{1}_n\}$ and therefore, as long as $\mathcal{G}$ is connected, the optimization problem above can be equivalently written as
\begin{align}
\label{eq:centralized_dual}
&\underset{\mathbf{y} \in \mathbb{R}^{np}}{\text{min}}  \quad f^*(\mathbf{y}) \\
&\text{s.t.} \quad  \left(\boldsymbol{\Pi}_{\mathbf{1}_n} \otimes \mathbf{I}_p \right) \mathbf{y} = \mathbf{0}, \nonumber
\end{align}
where $\boldsymbol{\Pi}_{\mathbf{1}_n} \triangleq \frac{1}{n}\mathbf{1}_n\mathbf{1}_n^T.$
Problem \eqref{eq:centralized_dual} can be solved using projected gradient descent yielding the iterates
\begin{subequations}
\begin{align}
\mathbf{x}(k+1) &:= \text{arg }\underset{\mathbf{x} \in \mathbb{R}^{np}}{\min} \quad f(\mathbf{x}) - \mathbf{y}(k)^T\mathbf{x} \\
\mathbf{y}(k+1) &:= \left(\boldsymbol{\Pi}^{\perp}_{\mathbf{1}_n} \otimes \mathbf{I}_p\right)\left(\mathbf{y}(k) - c \mathbf{x}(k+1)\right) \label{eq:grad_denscent_y},
\end{align}
\end{subequations}
where $\boldsymbol{\Pi}_{\mathbf{1}_n}^{\perp} \triangleq \mathbf{I}_n - \boldsymbol{\Pi}_{\mathbf{1}_n}.$ Note that if $f^*$ is differentiable $\mathbf{x}(k+1) = \nabla f^*(\mathbf{y}(k))$ \cite{rockafellar}. Further, if we select $\mathbf{y}(0) = (\boldsymbol{\Pi}_{\mathbf{1}_n}^{\perp} \otimes \mathbf{I}_p)\mathbf{y}(0)$ \eqref{eq:grad_denscent_y} can be rewritten as
\begin{align}
\mathbf{y}(k+1) := & \mathbf{y}(k) - c \left( \boldsymbol{\Pi}_{\mathbf{1}_n}^{\perp} \otimes \mathbf{I}_p\right) \mathbf{x}(k+1), \\
& =\mathbf{y}(k) - c \left( \mathbf{x}(k+1) - \left(\boldsymbol{\Pi}_{\mathbf{1}_n}\otimes \mathbf{I}_p\mathbf{x}(k+1)\right)\right)
\end{align}
which corresponds to using dual decomposition applied to the problem
\begin{align}
\label{eq:centralized}
& \underset{\mathbf{x} \in \mathbb{R}^{np}}{\text{min}} \quad f(\mathbf{x}) \\
& \text{s.t.} \quad \left( \boldsymbol{\Pi}^{\perp}_{\mathbf{1}_n} \otimes \mathbf{I}_p \right)\mathbf{x} = \mathbf{0}. \nonumber
\end{align}
Note that the optimization problem \eqref{eq:centralized} implicitly assumes that the graph $\mathcal{G}$ is fully connected. While the formulation \eqref{eq:centralized} is not practically useful, since it assumes all nodes can communicate with each other in a single hop, it illustrates that an equivalence with the fully connected graph can be drawn to any connected graph. This allows us to think of the formulation in \eqref{eq:centralized} as ``graph independent."  While using dual decomposition on \eqref{eq:centralized} yields iterates that are not computable in a distributed fashion unless the graph is fully connected, we can easily identify where the problematic term is. 

Assuming that $\mathcal{G}$ is not fully connected the only computation requiring communication between agents that are too far apart is $(\boldsymbol{\Pi}_{\mathbf{1}_n} \otimes \mathbf{I}_p) \mathbf{x}(k+1).$ Hence, we require a way to, in a distributed manner, compute or approximate $(\boldsymbol{\Pi}_{\mathbf{1}_n} \otimes \mathbf{I}_p) \mathbf{x}(k+1).$ Let $\{\mathbf{z}(k)\}_{k \geq 0}$ denote a sequence tasked with approximating $(\boldsymbol{\Pi}_{\mathbf{1}_n} \otimes \mathbf{I}_p)\mathbf{x}(k+1).$ With this, we have described the following method
\begin{subequations}
\begin{align}
\mathbf{x}(k+1) &:= \nabla f^*(\mathbf{y}(k)) \\
\mathbf{z}(k+1) &:= \,? \\
\mathbf{y}(k+1) &:= \mathbf{y}(k) - c \left(\mathbf{x}(k+1) - \mathbf{z}(k+1)\right), 
\end{align}
\end{subequations}
where we need a way of computing $\mathbf{z}(k+1)$ distributedly. 
To this end, let $\mathcal{Y} \triangleq \{\mathbf{y} \in \mathbb{R}^{np} : (\boldsymbol{\Pi}_{\mathbf{1}_n} \otimes \mathbf{I}_p)\mathbf{y} = \mathbf{0}\}.$ Further, let $[\cdot]_{\mathcal{Y}}$ denote the projection operator on the set $\mathcal{Y}.$  Before introducing which technique we use to approximate $(\boldsymbol{\Pi}_{\mathbf{1}_n} \otimes \mathbf{I}_p)\mathbf{x}(k+1),$ we will discuss what properties we want the iterate in $\mathbf{z}(k+1)$ to fulfil. First and foremost, we require that
\begin{equation}
\mathbf{z}(\infty) = \left(\boldsymbol{\Pi}_{\mathbf{1}_n} \otimes \mathbf{I}_p \right)\mathbf{x}(\infty).
\end{equation} 
Additional to $\mathbf{z}(\infty) = (\boldsymbol{\Pi}_{\mathbf{1}_n} \otimes \mathbf{I}_p)\mathbf{x}(\infty)$ we wish that as 
\begin{equation}
\|\mathbf{x}(k+1) - \mathbf{x}(k)\| \to 0,
\end{equation} 
it should also follow that
\begin{equation}\|\mathbf{z}(k+1) - (\boldsymbol{\Pi}_{\mathbf{1}_n} \otimes \mathbf{I}_p) \mathbf{z}(k+1)\| \to  0.
\end{equation} This is because, under certain assumptions, a bound on $\|\mathbf{x}(k+1) - \mathbf{x}(k)\|$ will imply a bound on $\|\mathbf{y}(k-1) - \left[\mathbf{y}(k-1) - c \mathbf{x}(k) \right]_{\mathcal{Y}}\|$ and vice-versa. Recall that the quantity 
\begin{equation}
\|\mathbf{y}(k-1) - \left[ \mathbf{y}(k-1) - c \mathbf{x}(k) \right]_{\mathcal{Y}}\|
\end{equation} is typically used to measure the sub-optimality gap for dual ascent \cite{luo}. Hence, what we are essentially requiring is that as we get close to optimality the approximation of $(\boldsymbol{\Pi}_{\mathbf{1}_n} \otimes \mathbf{I}_p)\mathbf{x}(k+1)$ should become better.

The work on DIGing in \cite{interpret} leverages the technique proposed in \cite{martinez} to distributedly solve \eqref{eq:original}. In the case of DIGing \cite{interpret}, the optimization problem in \eqref{eq:on_graph} is solved by averaging the primal variable, tracking the average gradient and performing steps in a descent direction. The average gradient tracking scheme in \cite{interpret} performs as follows. Let $\nabla f$ denote the gradient of the objective function in \eqref{eq:on_graph}. Then, the following scheme is used to track the quantity $(\boldsymbol{\Pi}_{\mathbf{1}_n}\otimes \mathbf{I}_p)\nabla f(\mathbf{x}(k+1))$ 
\begin{align}
\label{eq:gradient_tracking}
\mathbf{g}(k+1) &:= \left(\mathbf{W}(k) \otimes \mathbf{I}_p\right)\mathbf{g}(k)\\
 &+ \nabla f(\mathbf{x}(k+1)) - \nabla f(\mathbf{x}(k)), \nonumber
\end{align}
where $\mathbf{W}(k) \in \mathbb{R}^{n\times n}$ is a matrix defined on the edges of the graph $\mathcal{G}(k).$ The specific requirements on $\mathbf{W}(k)$ will be introduced in the following sub-section and are not necessary for an intuitive understanding of PANDA. Let $\mathbf{g}(k)$ be partitioned in $n$ sub-vectors of length $p,$ i.e. $\mathbf{g}(k) = [\mathbf{g}_1(k)^T,\hdots,\mathbf{g}_n(k)^T]^T.$ Then, the distributed nature of \eqref{eq:gradient_tracking} can be made obvious as
\begin{align}
\mathbf{g}_i(k+1) &:= \!\!\!\!\!\!\! \sum_{j \in \mathcal{N}_i(k) \cup \{i\}} \!\!\!\!\!\!\! w_{ij}(k)\mathbf{g}_j(k) \\&+ \nabla f_i(\mathbf{x}_i(k+1)) - \nabla f_i(\mathbf{x}_i(k)), \nonumber
\end{align} 
where $\nabla f_i$ denotes the gradient of the $i^{\text{th}}$ function in \eqref{eq:on_graph}.
As established in \cite{interpret} the scheme \eqref{eq:gradient_tracking} fulfills that $\mathbf{g}(\infty) = (\boldsymbol{\Pi}_{\mathbf{1}_n}\otimes \mathbf{I}_p)\nabla f(\infty).$ Further, as established in \cite{interpret,martinez} as 
\begin{equation}
\|\nabla f(\mathbf{x}(k+1)) - \nabla f(\mathbf{x}(k))\| \to 0,
\end{equation}
it follows that
\begin{equation}\|\mathbf{g}(k+1) - (\boldsymbol{\Pi}_{\mathbf{1}_n}\otimes \mathbf{I}_p)\mathbf{g}(k+1)\| \to 0.
\end{equation} Hence, the scheme in \eqref{eq:gradient_tracking} has all the properties we require. Therefore, in order to track the quantity $(\boldsymbol{\Pi}_{\mathbf{1}_n} \otimes \mathbf{I}_p)\mathbf{x}(k+1)$ we will use the scheme
\begin{equation}
\mathbf{z}(k+1) := \left(\mathbf{W}(k) \otimes \mathbf{I}_p\right)\mathbf{z}(k) + \mathbf{x}(k+1) - \mathbf{x}(k).
\end{equation}
We are with this ready to formally introduce PANDA which is written in a way that showcases its distributed nature in Algorithm \ref{alg:panda}.
\begin{algorithm}
\caption{PANDA \label{alg:panda}}
\begin{algorithmic}[1]
\State Choose step size $c > 0$ and pick $\mathbf{z}(0)= \mathbf{x}(0) = \mathbf{0}$ and $\mathbf{y}(0)$ such that $(\boldsymbol{\Pi}_{\mathbf{1}_n} \otimes \mathbf{I}_p)\mathbf{y}(0) = \mathbf{0}.$
\For{$k=0,1,\hdots$} each agent $i$:
\State computes $$
\mathbf{x}_i(k+1) := \text{arg } \underset{\mathbf{x}_i \in \mathbb{R}^{p}}{\text{min}}  f_i(\mathbf{x}_i) - \mathbf{y}_i(k)^T\mathbf{x}_i$$
\State exchanges $\mathbf{z}_i(k)$
 with the agents in $\mathcal{N}_i(k).$ 
\State computes 
$$\mathbf{z}_i(k+1) := \sum_{j \in \mathcal{N}_i(k) \cup \{i\}}w_{ij}(k)\mathbf{z}_j(k) + \mathbf{x}_i(k+1) - \mathbf{x}_i(k)$$
\State and computes $$\mathbf{y}_i(k+1) := \mathbf{y}_i(k) - c(\mathbf{x}_i(k+1) -\mathbf{z}_i(k+1))$$
\EndFor
\end{algorithmic}
\end{algorithm}

While the phrasing in Algorithm \ref{alg:panda} illustrates the algorithm's distributed nature it is not the most convenient way of representing it for analysis. We will therefore from now on use the following compact notation instead
\begin{subequations}
\begin{align}
\mathbf{x}(k+1) &:= \text{arg }\underset{\mathbf{x} \in \mathbb{R}^{np}}{\text{min}} \quad f(\mathbf{x}) - \mathbf{y}(k)^T\mathbf{x} \\
\mathbf{z}(k+1) &:= \left(\mathbf{W}(k)\otimes \mathbf{I}_p \right)\mathbf{z}(k) + \mathbf{x}(k+1) - \mathbf{x}(k) \label{eq:ziterate}\\
\mathbf{y}(k+1) &:= \mathbf{y}(k) - c (\mathbf{x}(k+1) - \mathbf{z}(k+1)). \label{eq:yiterate}
\end{align}
\end{subequations}

%%%%%%%%%%%%%%%%%%%%%%%%%%%%%%%%%%%%%%%%%%%%%%%%%%%%%%%%%%%%%%%%%%%%%%%%%%%
%%%%%%%%%%%%%%%%%%%%%%	PANDA R-LINEAR %%%%%%%%%%%%%%%%%%%%%%%%%%%%%%%%%%%%
%%%%%%%%%%%%%%%%%%%%%%%%%%%%%%%%%%%%%%%%%%%%%%%%%%%%%%%%%%%%%%%%%%%%%%%%%%
%%%%%%%%%%%%%%%%%%%%%%%%%%%%%%%%%%%
%%%%%%%%%% ASSUMPTION %%%%%%%%%%%%%
%%%%%%%%%%%%%%%%%%%%%%%%%%%%%%%%%%
\subsection{Convergence Statement}
In this sub-section we formally state the convergence properties of PANDA. However, before giving the formal statement we will provide the assumptions required for the statement to hold.
We first introduce some requirements on the sequence of mixing matrices $\{\mathbf{W}(k)\}.$
\begin{assumption}[Mixing matrix sequence $\{\mathbf{W}(k)\}$\cite{interpret}\label{assumption:mixing}]
For any $k=0,1,\hdots,$ the mixing matrix $\mathbf{W}(k) \in \mathbb{R}^{n \times n}$ satisfies the following relations:
\begin{enumerate}[label = (\subscript{P}{{\arabic*}})]
\item \emph{Decentralized property:} \label{prop:network} if $i \neq j$ and $(i,j) \not \in \mathcal{E}(k)$ $W_{ij}(k) = 0,$ i.e., $\mathbf{W}(k)$ is defined on the edges of the graph $\mathcal{G}(k).$
\item \emph{Double stochasticity:} \label{prop:consensus1} $\mathbf{W}(k)\mathbf{1}_n = \mathbf{1}_n,$ $\mathbf{1}_n^T\mathbf{W}(k) = \mathbf{1}_n^T.$
\item \emph{Joint spectrum property:} \label{prop:consensus2} Let 
\end{enumerate}
\begin{equation}
\mathbf{W}_b(k)\triangleq \mathbf{W}(k)\mathbf{W}(k-1)\hdots \mathbf{W}(k-b+1),
\end{equation}
for $k \geq 0$ and $b \geq k-1,$ with $\mathbf{W}_b(k) = \mathbf{I}_n$ for $k < 0$ and $\mathbf{W}_0(k) = \mathbf{I}_n.$ Then, there exists a positive $B$ such that
\begin{equation}
\underset{k \geq B-1}{\text{sup}} \quad \delta(k) = \delta < 1,
\end{equation}
where 
\begin{equation}
\delta(k) = \sigma_{\text{max}} \left\{ \mathbf{W}_B(k) - \frac{1}{n}\mathbf{1}_n^T\mathbf{1}_n \right\},\, \forall k=0,1,\hdots
\end{equation}
and $\sigma_{\text{max}}\left\{\cdot\right\}$ denotes the largest singular value of a matrix.
\end{assumption}
Properties \ref{prop:network} and \ref{prop:consensus1} in Assumption \ref{assumption:mixing} are common assumptions in the consensus literature  while \ref{prop:consensus2} is a requirement due to the connectivity (or lack of thereof) of the underlying graph $\mathcal{G}(k).$ Further discussion and examples under which Assumption \ref{assumption:mixing} holds can be found in \cite{interpret}, which also assumes \ref{prop:consensus2}.
We now introduce an assumption regarding the objective function in \eqref{eq:on_graph}.
%%%%%%%%%%%%%%%%%%%%%%%%%%%%%%%%
%%%%%%%%%%%ASSUMPTION %%%%%%%%%%
%%%%%%%%%%%%%%%%%%%%%%%%%%%%%%%%
\begin{assumption}[Strong convexity and smoothness\label{assumption:function}] The function $f$ is strongly convex and Lipschitz differentiable, i.e.,
\begin{equation}
f(\mathbf{x}) \geq f(\mathbf{y}) + (\nabla f(\mathbf{y}))^T(\mathbf{x} - \mathbf{y}) + \frac{\mu}{2}\|\mathbf{x} - \mathbf{y}\|^2
\end{equation}
and 
\begin{equation}
\|\nabla f(\mathbf{x}) - \nabla f(\mathbf{y})\| \leq L\|\mathbf{x} - \mathbf{y}\|,
\end{equation}
$\forall\, \mathbf{x},\,\mathbf{y} \in \mathbb{R}^{np}$ where $\mu > 0$ and $L < \infty$ are the strong convexity and Lipschitz constant respectively.
\end{assumption}
Assumption \ref{assumption:function} is a standard assumption to achieve geometric convergence. Some exceptions include \cite{non-strongly, polyak} (centralized) and \cite{extra, maros}(distributed).
%%%%%%%%%%%%%%%%%%%%%%%%%%%%%%
%%%%%%%%%% THEOREM %%%%%%%%%%%
%%%%%%%%%%%%%%%%%%%%%%%%%%%%%%
We are now ready to provide the paper's main statement regarding the convergence of PANDA. 
\begin{theorem}[PANDA converges R-linearly\label{theorem:panda}] Let Assumptions \ref{assumption:mixing} and \ref{assumption:function} hold. Also, let $\kappa \triangleq \frac{L}{\mu}$ denote the condition number of $f.$
Then for any step-size
\begin{equation}
c \in \left( 0 , \frac{\mu(1-\delta)^2}{2\sqrt{\kappa}}\right),
\end{equation}
the sequence $\{\mathbf{y}(k)\}$ converges to $\mathbf{y}^{\star},$ the unique solution of \eqref{eq:centralized_dual}, and $\{\mathbf{x}(k)\}$ converges to $\mathbf{x}^{\star},$ the unique solution of \eqref{eq:centralized}, at a global R-linear rate $\mathcal{O}(\lambda^k),$ where $\lambda < 1$ is given by
\begin{equation}
\lambda = 
\begin{cases}
\sqrt[2B]{1-\frac{c}{2L}} & \text{if } c \in \left(0,\bar{c} \right] \\
\sqrt[B]{\delta + \sqrt{\frac{2c(\sqrt{\kappa})}{\mu}}} & \text{if } c \in \left(\bar{c}, \frac{\mu (1-\delta)^2}{2\sqrt{\kappa}}\right),
\end{cases}
\end{equation} 
where
\begin{equation}
\bar{c} = \frac{\mu}{2}\left(\frac{16\kappa^{\frac{3}{2}}-4\kappa(1-\delta^2)}{(1+4\kappa^{\frac{3}{2}})^2} \right).
\end{equation}
\end{theorem}

\section{Proof of Theorem \ref{theorem:panda} \label{section:proof}}
This section is devoted to proving Theorem \ref{theorem:panda}. The section is divided into three subsections. The first is dedicated to establishing the structure of the proof, inspired by that in \cite{interpret}. As in \cite{interpret} we use the small gain theorem to establish a circle of arrows. The second subsection is dedicated to establishing each of the arrows in the circle. Finally in the third subsection brings the use of the small gain theorem together with the established arrows in order to conclude that PANDA converges R-linearly. 
\subsection{Proof structure \label{subsection:proofstructure}}
We start by introducing some notation and the small gain theorem for completeness. Let $\mathbf{s}^i \triangleq \{\mathbf{s}^i(0),\mathbf{s}^i(1),\hdots\}$ denote an infinite sequence of vectors $\mathbf{s}^i(k) \in \mathbb{R}^{np},$ for $i = 1,\hdots m.$ Further, let
\begin{equation}
\|\mathbf{s}^i\|^{\lambda,K} \triangleq \underset{k=0,\hdots,K}{\text{sup}} \quad \frac{1}{\lambda^k}\|\mathbf{s}^i(k)\|
\end{equation}
and 
\begin{equation}
\|\mathbf{s}^i\|^{\lambda} \triangleq \underset{k \geq 0}{\text{sup}} \quad \frac{1}{\lambda^k}\|\mathbf{s}^i(k)\|.
\end{equation}
\begin{theorem}[small gain theorem \cite{interpret}\label{theorem:gain}]
Suppose $\mathbf{s}^1,\hdots,\mathbf{s}^m$ are vector sequences such that for all positive integers $K$ and for each $i = 1,\hdots,m,$ we have an arrow $\mathbf{s}^i \to \mathbf{s}^{(i\!\!\mod m) + 1},$ i.e.,
\begin{equation}
\label{eq:cyclic_gain}
\|\mathbf{s}^{(i \!\!\!\!\mod m) + 1}\|^{\lambda,K} \leq \gamma_i \|\mathbf{s}^i\|^{\lambda,K} + \omega_i,
\end{equation}
where the constants $\gamma_1,\hdots,\gamma_m$ and $\omega_1,\hdots,\omega_m$ are independent of $K.$ Further, suppose that the constants $\gamma_1,\hdots,\gamma_m$ are nonnegative and satisfy
\begin{equation}
\gamma_1\hdots\gamma_m < 1.
\end{equation}
Then we have that
\begin{align}
\|\mathbf{s}^1\|^{\lambda} \leq &\frac{1}{1 - \gamma_1\gamma_2 \hdots \gamma_m}\left(\omega_1 \gamma_2 \gamma_3 \hdots \gamma_m + \omega_2 \gamma_3\gamma_4 \hdots \gamma_m \right. \nonumber \\
& \left. + \hdots + \omega_{m-1}\gamma_m + \omega_m \right).  
\end{align}
\end{theorem}
Note that if $\|\mathbf{s}^1\|^{\lambda} \leq C$ where $C < \infty,$ $\|\mathbf{s}^1(k)\|$ converges to zero exponentially fast and at rate $\lambda.$ Proof of this statement can be found in \cite{interpret}. Further, due to the cyclic nature of \eqref{eq:cyclic_gain} all sequences $\{\|\mathbf{s}_i(k)\|\}$ converge to zero exponentially fast at rate $\lambda$ \cite{interpret}.

The goal is ultimately to establish that $\|\mathbf{y}(k) - \mathbf{y}^{\star}\| \to 0$ as $k \to \infty$ exponentially fast. For notational convenience let $\mathbf{r}(k) \triangleq \mathbf{y}(k) - \mathbf{y}^{\star}.$ In order to directly apply Theorem \ref{theorem:gain} we will start the circle of arrows with the sequence $\{\mathbf{r}(k)\}.$ We will then proceed as
\begin{equation}
\label{eq:cycle}
\mathbf{r} \to \mathbf{x}^{\perp} \to \boldsymbol{\Delta}_{xz}^{\perp} \to \Delta \mathbf{y} \to \mathbf{z}^{\perp} \to \mathbf{r},
\end{equation}
where 
\begin{subequations}
\begin{align}
\mathbf{r}(k) \triangleq & \,\,\mathbf{y}(k) - \mathbf{y}^{\star}, \, k \geq 0 \\
\mathbf{x}^{\perp}(k) \triangleq & \left( \boldsymbol{\Pi}_{\mathbf{1}_n}^{\perp}\otimes \mathbf{I}_p\right)\mathbf{x}(k) , \, k \geq 0 \label{eq:seqx}\\
\boldsymbol{z}^{\perp}(k) \triangleq & \left( \boldsymbol{\Pi}_{\mathbf{1}_n}^{\perp} \otimes \mathbf{I}_p \right)\mathbf{z}(k),\, k \geq 0 \\
\Delta\mathbf{y}(k) \triangleq & \,\,\mathbf{y}(k) - \mathbf{y}(k-1),\, k \geq 0 \\
\boldsymbol{\Delta}_{xz}^{\perp}(k) \triangleq & \,\,\mathbf{x}^{\perp}(k) - \mathbf{z}^{\perp}(k), k \geq 0.
\end{align}
\end{subequations}
We adopt the convention that $\mathbf{x}(0) = \mathbf{z}(0) = \mathbf{0}$ and $\Delta \mathbf{y}(0) = \mathbf{0}$ which is consistent with the implementation in Algorithm 1.
More specifically, when it comes to the arrows \eqref{eq:cycle} we establish the following relations
\begin{itemize}
\item[(A1)] $\|\mathbf{x}^{\perp}\|^{\lambda,K} \leq \gamma_1 \|\mathbf{r}\|^{\lambda,K} + \omega_1,$ where $\gamma_1 = \frac{1}{\mu\lambda}$ $\omega_1 = 0.$
\item[(A2)] $\|\boldsymbol{\Delta}_{xz}^{\perp}\|^{\lambda,K} \leq \gamma_2 \|\mathbf{x}^{\perp}\|^{\lambda,K} + \omega_2,$ where $\gamma_2 = \frac{2(1-\lambda^B)}{(1-\lambda)(\lambda^B -\delta)}$ and $\omega_2 = \frac{\lambda^B}{\lambda^B - \delta} \sum_{t=1}^B \lambda^{1-t}\|\boldsymbol{\Delta}_{xz}^{\perp}(t-1)\|$ for $\lambda^B > \delta.$
\item[(A3)] $\|\Delta \mathbf{y}\|^{\lambda,K} \leq \gamma_3\|\boldsymbol{\Delta}_{xz}^{\perp}\|^{\lambda,K} + \omega_3,$ where $\gamma_3 = c$ and $\omega_3 = 0.$
\item[(A4)] $\|\mathbf{z}^{\perp}\|^{\lambda,K} \leq \gamma_4 \|\Delta \mathbf{y}\|^{\lambda,K} + \omega_4,$ where $\gamma_4 = \frac{(1-\lambda^B)}{\mu (1-\lambda)(\lambda^B - \delta)}$ and $\omega_4 = \frac{\lambda^B}{\lambda^B - \delta} \sum_{t=1}^B \lambda^{1-t}\|\mathbf{z}(t-1)\|$ for $\lambda^B > \delta.$
\item[(A5)] $\|\mathbf{r}\|^{\lambda,K} \leq \gamma_5\|\mathbf{z}^{\perp}\|^{\lambda,K} + \omega_5,$ where $\gamma_5 = \sqrt{L\mu}$ and $\omega_5 = 2 \|\mathbf{r}(0)\|$ for $\lambda \in [\sqrt{1-\frac{c}{2L}},1)$ and $c \in (0,\frac{\mu}{2}].$
\end{itemize}
Once all of these relations are established it remains to show that there exists pairs $(\lambda,c)$ such that $\lambda < 1,$ $\gamma_1\hdots\gamma_5 < 1$ and all the restrictions regarding $\lambda$ and $c$ included in (A1)-(A5) hold.

The next subsection is devoted to establishing the circle of arrows in \eqref{eq:cycle}. This will be done in the form of Lemmas \ref{lemma:a1}-\ref{lemma:a4}, and \ref{lemma:a5}. However, before proceeding to establishing the arrows above we require some additional results. These results are stated below and relate to the function $f^*,$ the averaging effects of the sequence of matrices $\{\mathbf{W}(k)\}$ and PANDA's iterates.
%%%%%%%%%%%%%%%%%%%%%%%%%%%%%%%%%%%%%%%%%%%%%%%%%%%%%%%%%%
%%%%%%%%%%%%%%%%%% THEOREM DUAL %%%%%%%%%%%%%%%%%%%%%%%%%%
%%%%%%%%%%%%%%%%%%%%%%%%%%%%%%%%%%%%%%%%%%%%%%%%%%%%%%%%%%
\begin{theorem}[Dual smoothness and strong convexity\label{theorem:dual}]
The function $f$ is strongly convex with constant $\mu$ if and only if its convex conjugate $f^*$ has Lipschitz gradients with constant $\frac{1}{\mu}.$ Further, the function $f$ has Lipschitz gradients with constant $L$ if and only if the function $f^*$ is strongly convex with constant $\frac{1}{L}.$
\end{theorem}
\begin{proof}
The proof follows by combining the proof of Theorem 6 in \cite{toyota} and Theorem 2.1.5. in \cite{nesterov}.
\end{proof}
%%%%%%%%%%%%%%%%%%%%%%%%%%%%%%%%%%%%%%%%%%%%%%%%%%%%%%%%%
%%%%%%%%%%%%%%%%%% LEMMA MIXING %%%%%%%%%%%%%%%%%%%%%%%%%
%%%%%%%%%%%%%%%%%%%%%%%%%%%%%%%%%%%%%%%%%%%%%%%%%%%%%%%%%
\begin{lemma}[$B-$step consensus contraction\label{lemma:mixing}]
Under Assumption \ref{assumption:mixing}, for any $k= B-1,B,\hdots,$ and any vector $\mathbf{b} \in \mathbb{R}^{np}$ if $\mathbf{a} = (\mathbf{W}_{B}(k) \otimes \mathbf{I}_p)\mathbf{b},$ then
\begin{equation}
\left\|(\boldsymbol{\Pi}_{\mathbf{1}_n}^{\perp} \otimes \mathbf{I}_p)\mathbf{a}\right\| \leq \delta(k)\left\|(\boldsymbol{\Pi}_{\mathbf{1}_n}^{\perp} \otimes \mathbf{I}_p)\mathbf{b}\right\|.
\end{equation}
\end{lemma}
The proof of lemma follows from Assumption \ref{assumption:mixing}, and a complete proof of the lemma can be found in \cite{interpret}.
\begin{lemma}[Equivalent iterates\label{lemma:equivalent}]
Under Assumption \ref{assumption:mixing} the PANDA iterate \eqref{eq:yiterate} can be equivalently written as 
\begin{equation}
\mathbf{y}(k+1) := \mathbf{y}(k) - c (\boldsymbol{\Pi}_{\mathbf{1}_n}^{\perp} \otimes \mathbf{I}_p)(\mathbf{x}(k+1) - \mathbf{z}(k+1)).
\end{equation} 
\end{lemma}
\begin{proof}
In order to see this, let us express $\mathbf{x}(k+1) - \mathbf{z}(k+1)$ exclusively as  a function of the iterates $\{\mathbf{x}(t)\}_{t=0}^k.$ Recall that the variables are initialized as $\mathbf{x}(0) = \mathbf{z}(0).$ Consider the PANDA iterate in \eqref{eq:ziterate} and rewrite it as
\begin{align}
\label{eq:l2it}
\mathbf{x}(k+1) - \mathbf{z}(k+1) = \left( \left( \mathbf{I}_n - \mathbf{W}(k) \right) \otimes \mathbf{I}_p\right)\mathbf{x}(k)  \\
 + \left( \mathbf{W}(k) \otimes \mathbf{I}_p \right)(\mathbf{x}(k) - \mathbf{z}(k)), \nonumber
\end{align}
which can be done by re-arranging \eqref{eq:ziterate} and adding and subtracting the quantity $(\mathbf{W}(k) \otimes \mathbf{I}_p)\mathbf{x}(k).$ The analogous procedure can be applied to the quantity $\mathbf{x}(k) - \mathbf{z}(k),$ i.e.,
\begin{align}
\mathbf{x}(k) - \mathbf{z}(k) = \left( \left( \mathbf{I}_n - \mathbf{W}(k-1) \right) \otimes \mathbf{I}_p\right)\mathbf{x}(k-1)  \\
 + \left( \mathbf{W}(k-1) \otimes \mathbf{I}_p \right)(\mathbf{x}(k-1) - \mathbf{z}(k-1)), \nonumber
\end{align}
which we evaluate in \eqref{eq:l2it} yielding
\begin{align}
\mathbf{x}(k+1) - \mathbf{z}(k+1) = \left( \left(\mathbf{I}_n - \mathbf{W}(k)\right) \otimes \mathbf{I}_p \right)\mathbf{x}(k) \\
+ \left( \mathbf{W}(k)\left( \mathbf{I}_n - \mathbf{W}(k-1) \right) \otimes \mathbf{I}_p\right)\mathbf{x}(k-1) \nonumber\\
+ (\mathbf{W}_{2}(k) \otimes \mathbf{I}_p)(\mathbf{x}(k-1) - \mathbf{z}(k-1)). \nonumber
\end{align}
By re-writing $\mathbf{x}(k-1)-\mathbf{z}(k-1),\hdots,\mathbf{x}(0) - \mathbf{z}(0)$ as they appear in the expression we obtain
\begin{align}
\mathbf{x}(k+1) & -  \mathbf{z}(k+1) =  \left(\left(\mathbf{I}_n - \mathbf{W}(k)\right)\otimes \mathbf{I}_p\right) \mathbf{x}(k)  \\
& +\left(\mathbf{W}(k)\left(\mathbf{I}_n - \mathbf{W}(k-1)\right)\otimes \mathbf{I}_p\right)\times\mathbf{x}(k-1) \nonumber\\ & + \sum_{t=0}^{k-3}\left(\mathbf{W}_{k-1-t}(k) \nonumber \left(\mathbf{I}_n - \mathbf{W}(t+1)\right)\otimes \mathbf{I}_p\right)\mathbf{x}(t+1) \nonumber,
\end{align}
where we use the fact that $\mathbf{x}(0) - \mathbf{z}(0) = \mathbf{0}.$
Since the matrices in the sequence $\{\mathbf{W}(k)\}$ are column stochastic it follows that $ \boldsymbol{\Pi}_{\mathbf{1}_n}^{\perp}\mathbf{W}_b(k)(\mathbf{I}_n - \mathbf{W}(t)) = \mathbf{W}_b(k)(\mathbf{I}_n - \mathbf{W}(t))$ for any $b \geq 0$ and $t \geq 0.$
\end{proof}

\subsection{Establishing the circle of arrows}
Before diving into the proof technicalities we will provide some additional intuition. Intuitively speaking, most convex optimization algorithms establish Q-linear convergence, under Assumption \ref{assumption:function}, by using, directly or indirectly, that the optimal point $\mathbf{x}^{\star}$ can be seen as the unique fixed point of a strictly contractive mapping. This is due to the fact that many algorithms can be seen as fixed point methods.

 In order to deal with time-varying networks, both DIGing and PANDA introduce terms that are aimed at tracking the average of a quantity. Since this can not be done accurately at all iterations due to the decentralized nature of the problem, each iteration inherently carries with them a certain amount of error. Hence, instead of treating the iterates purely as fixed point iterates, both DIGing and PANDA need to make the claim that as the algorithm progresses, the approximations become better. In particular, the iterates of PANDA can be seen as iterates of projected gradient descent on \eqref{eq:centralized_dual} where the gradient is erroneously computed (c.f. Lemma \ref{lemma:equivalent}).
   The circle of arrows as formulated in \eqref{eq:cycle} reflects this intuition. More specifically, we start by assuming that $\|\mathbf{r}(k)\| \to 0$ exponentially at rate $\lambda.$ In turn, this has several implications, among which are that the error incurred at each iteration of PANDA is a geometrically decreasing quantity. Once this has been established, the circle is closed by using similar procedures as other standard convex optimization methods, with one essential difference. As will be seen in the proof of Lemma \ref{lemma:a5}, PANDA is most sensitive to large errors as it approaches the optimal point. However, since the error vanishes as the algorithm progresses this does not jeopardize the algorithm's convergence as long as the step-size is chosen to be sufficiently small.  
How small the steps need to be is established in subsection \ref{subsection:small} by selecting pairs $(\lambda,c)$ that fulfill the requirements of the small gain theorem. Intuitively speaking, additional to the standard requirements on step-size, the step-size must be chosen small enough so that the algorithm progresses slowly enough for the tracking scheme to be able to keep up. With this said, we now proceed to establishing the circle of arrows.

%%%%%%%%%%%%%%%%%%%%%%%%%%%%%%%%%%%%%%%%%%%%%%%%%%%%%%%%%%
%%%%%%%%%%%%%%%%%% ARROW 1 %%%%%%%%%%%%%%%%%%%%%%%%%%%%%%%
%%%%%%%%%%%%%%%%%%%%%%%%%%%%%%%%%%%%%%%%%%%%%%%%%%%%%%%%%%
\begin{lemma}[A1\label{lemma:a1}]
If $f$ is strongly convex with constant $\mu$ (Assumption \ref{assumption:function}), it holds that
\begin{equation}
\|\mathbf{x}^{\perp}\|^{\lambda,K} \leq \|\mathbf{r}\|^{\lambda,K},
\end{equation}
for all $\lambda \in (0,1)$ and $K \geq 0.$
\end{lemma}
\begin{proof}
From Theorem 3 it follows that $f^*$ has Lipschitz gradients with constant $\frac{1}{\mu}.$ Further, since $x^{\star} = \nabla f^{*}(\mathbf{y}^{\star})$ and $\mathbf{x}(k+1) = \nabla f^*(\mathbf{y}(k))$ it follows that
\begin{equation}
\|\mathbf{x}(k+1) - \mathbf{x}^{\star}\| \leq \frac{1}{\mu}\|\mathbf{y}(k) - \mathbf{y}^{\star}\|.
\end{equation}
By taking the square of $\|\mathbf{x}(k+1) - \mathbf{x}^{\star}\|$ and using the fact that $\mathbf{x}(k+1) = (\boldsymbol{\Pi}_{\mathbf{1}_n} \otimes \mathbf{I}_p)\mathbf{x}(k+1)  + (\boldsymbol{\Pi}_{\mathbf{1}_n}^{\perp} \otimes \mathbf{I}_p)\mathbf{x}(k+1)$ we obtain
\begin{align}
\left\|\mathbf{x}(k+1) - \mathbf{x}^{\star} \right\|^2  =  \left\|\boldsymbol{\Pi}_{\mathbf{1}_n}^{\perp} \mathbf{x}(k+1)\right\|^2   + \left\|\boldsymbol{\Pi}_{\mathbf{1}_n}\mathbf{x}(k+1)\right\|^2
\end{align}
where we have used that $((\boldsymbol{\Pi}_{\mathbf{1}_n}^{\perp} \otimes \mathbf{I}_p)\mathbf{x}(k+1))^T((\boldsymbol{\Pi}_{\mathbf{1}_n} \otimes \mathbf{I}_p)\mathbf{x}(k+1) - \mathbf{x}^{\star}) = 0.$ Consequently, it holds that
\begin{equation}
\left\| \left(\boldsymbol{\Pi}_{\mathbf{1}_n}^{\perp} \otimes \mathbf{I}_p \right)\mathbf{x}(k+1)\right\| \leq \|\mathbf{x}(k+1)-\mathbf{x}^{\star} \|
\end{equation}
which implies (c.f. \eqref{eq:seqx})
\begin{equation}
\label{eq:bothsides}
\lambda^{-k}\|\mathbf{x}^{\perp}(k)\| \leq  \frac{1}{\mu\lambda}\lambda^{-k}\|\mathbf{r}(k)\|.
\end{equation}
Note that
\begin{align}
\|\mathbf{x}^{\perp}\|^{\lambda,K} = & \underset{k=-1,\hdots,K-1}{\sup} \quad \lambda^{-(k+1)}\|\mathbf{x}^{\perp}(k+1)\| \nonumber\\
& = \underset{k=0,\hdots,K-1}{\sup} \quad \lambda^{-(k+1)}\|\mathbf{x}^{\perp}(k+1)\|,
\end{align}
because $\|\mathbf{x}^{\perp}(0)\|=0.$
Hence, taking $\!\!\!\!\!\!\underset{k=0,1,\hdots,K}{\text{sup}}\!\!\!\!\!\!$ on both sides of \eqref{eq:bothsides} yields
\begin{equation}
\|\mathbf{x}^{\perp}\|^{\lambda,K} \leq \|\frac{1}{\lambda \mu}\mathbf{r}\|^{\lambda,K-1} \leq \frac{1}{\mu \lambda}\|\mathbf{r}\|^{\lambda,K}.
\end{equation}
\end{proof}
%%%%%%%%%%%%%%%%%%%%%%%%%%%%%%%%%%%%%%%%%%%%%%%%%%%%%%%%%%%%%%%
%%%%%%%%%%%%%%%%%%%%%%%% ARROW 2 %%%%%%%%%%%%%%%%%%%%%%%%%%%%%%
%%%%%%%%%%%%%%%%%%%%%%%%%%%%%%%%%%%%%%%%%%%%%%%%%%%%%%%%%%%%%%%
\begin{lemma}[A2\label{lemma:a2}]
 Under Assumption \ref{assumption:mixing} it holds that
\begin{align}
\|\boldsymbol{\Delta}_{xz}^{\perp}\|^{\lambda,K} \leq & \frac{2(1-\lambda^B)}{(1-\lambda)(\lambda^B - \delta)}\|\mathbf{x}^{\perp}\|^{\lambda,K}  \\
& +\frac{\lambda^B}{\lambda^B - \delta}\sum_{t=1}^B \lambda^{1-t}\|\boldsymbol{\Delta}_{xz}^{\perp}(t-1)\|.
\end{align}
for all $\lambda \in (\delta^{\frac{1}{B}},1)$ and for all $K \geq 0.$
\end{lemma}

\begin{proof}
By using that $\mathbf{z}(k+1) = (\mathbf{W}(k)\otimes\mathbf{I}_p)\mathbf{z}(k) + \mathbf{x}(k+1) - \mathbf{x}(k)$ (c.f. \eqref{eq:ziterate}) and adding and subtracting $(\boldsymbol{\Pi}_{\mathbf{1}_n}^{\perp}\mathbf{W}(k) \otimes \mathbf{I}_p)\mathbf{x}(k),$ we have
\begin{align}
\label{eq:iterations}
& \left(\boldsymbol{\Pi}_{\mathbf{1}_n}^{\perp} \otimes \mathbf{I}_p\right)\left(\mathbf{z}(k+1) - \mathbf{x}(k+1) \right) = \left(\mathbf{W}(k)\boldsymbol{\Pi}_{\mathbf{1}_n}^{\perp}\otimes \mathbf{I}_p\right) \nonumber \\
&\times \left(\mathbf{z}(k) - \mathbf{x}(k)\right)-\left(\boldsymbol{\Pi}_{\mathbf{1}_n}^{\perp}\left(\mathbf{I}_n - \mathbf{W}(k)\right) \otimes \mathbf{I}_p\right)\mathbf{x}(k).
\end{align}
Note that we have used that the matrices $\boldsymbol{\Pi}_{\mathbf{1}_n}^{\perp}$ and $\mathbf{W}(k)$ commute because $\mathbf{W}(k)$ is doubly stochastic. Note that on the RHS of \eqref{eq:iterations} we have the quantity $(\boldsymbol{\Pi}_{\mathbf{1}_n}^{\perp}\otimes \mathbf{I}_p)(\mathbf{z}(k) - \mathbf{x}(k)).$ Due to the iterative nature of PANDA we can express  $(\boldsymbol{\Pi}_{\mathbf{1}_n}^{\perp}\otimes \mathbf{I}_p)(\mathbf{z}(k) - \mathbf{x}(k))$ as \eqref{eq:iterations} by simply shifting the time indexes one step backward. Then, the obtained expression can be substituted back into \eqref{eq:iterations}. By iteratively applying this procedure $B$ times, we obtain that for any $k \geq B-1$ we have that
\begin{align}
&\left(\boldsymbol{\Pi}_{\mathbf{1}_n}^{\perp} \otimes \mathbf{I}_p\right)\left(\mathbf{z}(k+1) - \mathbf{x}(k+1)\right) = 
\left(\boldsymbol{\Pi}_{\mathbf{1}_n}^{\perp}\mathbf{W}_B(k) \otimes \mathbf{I}_p\right)\nonumber \\
&\quad \times\left(\mathbf{z}(k-(B-1))- \mathbf{x}(k-(B-1))\right) \nonumber \\
& \quad - \sum_{t=0}^{B-1}(\boldsymbol{\Pi}_{\mathbf{1}_n}^{\perp}\mathbf{W}_{t}(k)(\mathbf{I}_n - \mathbf{W}(k-t))\otimes \mathbf{I}_p)\mathbf{x}(k-t).
\end{align}
Let $\|\cdot\|_{\boldsymbol{\Pi}^{\perp}_{\mathbf{1}_n} \otimes \mathbf{I}_p}$ be the semi-norm induced by the matrix $\boldsymbol{\Pi}^{\perp}_{\mathbf{1}_n} \otimes \mathbf{I}_p.$ For convenience, let $l \triangleq k-(B-1).$ Then
\begin{align}
\label{eq:long_recursion}
&\left\|\left(\boldsymbol{\Pi}_{\mathbf{1}_n}^{\perp} \otimes \mathbf{I}_p\right)\left(\mathbf{z}(k+1) - \mathbf{x}(k+1)\right) \right\|_{\boldsymbol{\Pi}_{\mathbf{1}_n}^{\perp} \otimes \mathbf{I}_p}  \leq \\  
&\qquad\left\|\left(\mathbf{W}_B(k)\boldsymbol{\Pi}_{\mathbf{1}_n}^{\perp} \otimes \mathbf{I}_p\right)(\mathbf{z}(l)  - \mathbf{x}(l))\right\|_{\boldsymbol{\Pi}_{\mathbf{1}_n}^{\perp} \otimes \mathbf{I}_p} \nonumber \\
& \!\!\!\!+ \sum_{t=0}^{B-1}\left\|\left(\mathbf{W}_t(k)(\mathbf{I}_n - \mathbf{W}(k-t))\boldsymbol{\Pi}_{\mathbf{1}_n}^{\perp} \otimes \mathbf{I}_p\right)\mathbf{x}(k-t)\right\|_{\boldsymbol{\Pi}_{\mathbf{1}_n}^{\perp} \otimes \mathbf{I}_p} \nonumber
\end{align}
where we have used that the matrices $\mathbf{W}_t(k)$ commute with  $\boldsymbol{\Pi}_{\mathbf{1}_n}^{\perp}$ for all $k$ and $t,$ $(\mathbf{I}_n - \mathbf{W}(k-t))$ commutes with $\boldsymbol{\Pi}_{\mathbf{1}_n}^{\perp}$ for all $k$ and $t,$ and the matrix $\boldsymbol{\Pi}_{\mathbf{1}_n}^{\perp}$ is idempotent. We can now rewrite \eqref{eq:long_recursion} as
\begin{align}
\label{eq:until_here}
\|\boldsymbol{\Delta}_{xz}^{\perp}(k+1)\| \leq 
 \left\|\left(\boldsymbol{\Pi}_{\mathbf{1}_n}^{\perp}\mathbf{W}_B(k) \otimes \mathbf{I}_p \right) \boldsymbol{\Delta}_{xz}^{\perp}(k-(B-1))\right\| \nonumber \\
 \sum_{t=0}^{B-1}\left\| \left( \mathbf{W}_t(k)(\mathbf{I}_n - \mathbf{W}(k-t))\boldsymbol{\Pi}_{\mathbf{1}_n}^{\perp} \otimes \mathbf{I}_p \right)\mathbf{x}(k-t) \right\|,
\end{align}
where we have used that $\rho(\boldsymbol{\Pi}_{\mathbf{1}_n}^{\perp}) = 1,$ where $\rho(\cdot)$ denotes a matrix's spectral radius. Since the matrices $\{\mathbf{W}(k)\}$ are doubly stochastic we have that $\rho(\mathbf{W}_t(k))\rho(\mathbf{I}_n - \mathbf{W}(k-t)) \leq 2.$ By further using that $\boldsymbol{\Pi}_{\mathbf{1}_n}^{\perp}\mathbf{W}_B(k) = \mathbf{W}_B(k)\boldsymbol{\Pi}_{\mathbf{1}_n}^{\perp}$ due to the double stochasticity of $\{\mathbf{W}_B(k)\}$ we can use Lemma \ref{lemma:mixing} to establish
\begin{align}
\label{eq:afterlemma}
\|\boldsymbol{\Delta}_{xz}^{\perp}(k+1)\| \leq 
& \delta \left\|\boldsymbol{\Delta}_{xz}^{\perp}(k-(B-1))\right\| \\
& + \sum_{t=1}^{B}2\left\|\mathbf{x}^{\perp}(k + 1-t) \right\|, \quad k \geq B-1. \nonumber
\end{align}
From here the procedure is identical to that in the proof of Lemma 3.10 in \cite{interpret}. We include the proof here for completeness. For $k = -1,\hdots, B-2$ we have that
\begin{align}
\label{eq:inequality_1}
\underset{k=-1,\hdots,B-2}{\text{sup}}  \lambda^{-(k+1)}\|\boldsymbol{\Delta}_{xz}^{\perp}(k+1)\| \leq \sum_{t=1}^B  \lambda^{1-t}\|\boldsymbol{\Delta}_{xz}^{\perp}(1-t)\|.
\end{align}
Note that this follows because the supremum of a finite positive sequence needs to be upper bounded by the sum of its elements.
Further, by taking the supremum on both sides of \eqref{eq:afterlemma} for $k \geq B-1$ we have
\begin{align}
\label{eq:inequality_2}
&\underset{k = B-1,\hdots,K-1}{\text{sup}} \lambda^{-(k+1)}\|\boldsymbol{\Delta}_{xz}^{\perp}(k+1)\| \leq \\
& \,\, \underset{k= B-1,\hdots,K-1}{\text{sup}} \frac{\delta}{\lambda^{B}}\lambda^{-(k+1-B)}   \|\boldsymbol{\Delta}_{xz}^{\perp}(k+1-B)\| \nonumber \\
& \,\,+ \underset{k=B-1,\hdots,K-1}{\text{sup}}\lambda^{-(k+1)}\sum_{t=1}^B2\|\mathbf{x}^{\perp}(k+1-t)\|. \nonumber
\end{align}
Also, since $\|\boldsymbol{\Delta}_{xz}^{\perp}\|^{\lambda,K} = \underset{k = -1,\hdots,K-1}{\text{sup}} \|\boldsymbol{\Delta}_{xz}^{\perp}(k+1)\| $ we have
\begin{align}
\label{eq:60}
\|\boldsymbol{\Delta}_{xz}^{\perp}\|^{\lambda,K} \leq &  
\underset{k = -1,\hdots, B-2}{\text{sup}} \, \lambda^{-(k+1)}\|\boldsymbol{\Delta}_{xz}^{\perp}(k+1)\| \\ &+ \!\!\!\! \underset{k = B-1,\hdots, K-1}{\text{sup}} \lambda^{-(k+1)}\|\boldsymbol{\Delta}_{xz}^{\perp}(k+1)\|. \nonumber
\end{align}
We can replace the quantities in the RHS of \eqref{eq:60} above  by their upper bounds using \eqref{eq:inequality_1} and \eqref{eq:inequality_2} yielding
\begin{align}
&\|\boldsymbol{\Delta}_{xz}^{\perp}\|^{\lambda,K} \leq \frac{\delta}{\lambda^B} \underset{k = -1,\hdots,B-1}{\text{sup}}\lambda^{-(k+1-B)}\|\boldsymbol{\Delta}^{\perp}(k+1-B)\| \nonumber \\
& \,\,\,+ \lambda^{-t} \underset{k=-1,\hdots,B-1}{\text{sup}} \lambda^{-(k+1-t)} \sum_{t=1}^B 2 \|\mathbf{x}^{\perp}(k+1-t)\| \\
&\,\,\,+ \sum_{t=1}^B \lambda^{1-t}\|\boldsymbol{\Delta}_{xz}^{\perp}(1-t)\|,\nonumber
\end{align}
which leads to
\begin{align}
\|\boldsymbol{\Delta}_{xz}^{\perp}\|^{\lambda,K} \leq &  \frac{\delta}{\lambda^B}\|\boldsymbol{\Delta}_{xz}^{\perp}\|^{\lambda,K} + 2 \sum_{t=1}^B \lambda^{-t}\|\mathbf{x}^{\perp}\|^{\lambda,K} + \\ & \sum_{t=1}^B \lambda^{1-t}\|\boldsymbol{\Delta}_{xz}^{\perp}(1-t)\| \nonumber,
\end{align}
which after some algebraic manipulations yields
\begin{align}
\label{eq:factor2}
\|\boldsymbol{\Delta}_{xz}^{\perp}\|^{\lambda,K} \leq &  \frac{2(1-\lambda^B)}{(1-\lambda)(\lambda^B - \delta)}\|\mathbf{x}^{\perp}\|^{\lambda,K} \\ & + \frac{\lambda^B}{\lambda^B - \delta}\sum_{t=1}^B \lambda^{1-t}\|\boldsymbol{\Delta}_{xz}^{\perp}(1-t)\| \nonumber,
\end{align}
establishing the desired result.
\end{proof}
Note that if the sequence of matrices $\{\mathbf{W}(k)\}$ consists of positive semi-definite matrices, i.e. $\mathbf{W}(k) \succeq \mathbf{0}$ then $\rho(\mathbf{W}_t(k)(\mathbf{I}_n- \mathbf{W}(k-t))) \leq \rho(\mathbf{W}_t(k))\rho(\mathbf{I}_n - \mathbf{W}(k-t)) \leq 1$ and the factor 2 in \eqref{eq:factor2} can be removed.
\begin{lemma}[A3] It holds that
\begin{equation}
\|\Delta \mathbf{y}\|^{\lambda,K} \leq c \|\boldsymbol{\Delta}_{xz}^{\perp}\|^{\lambda,K}.
\end{equation}
for any $c > 0,$ $\lambda \in (0,1)$ and $K \geq 0.$
\end{lemma}
\begin{proof}
The proof follows easily by using that $\mathbf{y}(k+1) - \mathbf{y}(k) = -c \boldsymbol{\Delta}_{xz}^{\perp}(k+1)$. Note that the bound is also fulfilled for $K=0$ because all involved sequences are initialized at $\mathbf{0}$.
\end{proof}
\begin{lemma}[A4\label{lemma:a4}]
Under Assumptions \ref{assumption:mixing} and \ref{assumption:function} it holds that
\begin{align}
\|\mathbf{z}^{\perp}\|^{\lambda,K} \leq \frac{(1-\lambda^B)}{\mu(1-\lambda)(\lambda^B - \delta)}\|\Delta \mathbf{y}\|^{\lambda,K} + \\
\frac{\lambda^B}{\lambda^B - \delta}\sum_{t=1}^B \lambda^{1-t}\|\mathbf{z}^{\perp}(t-1)\|. \nonumber
\end{align}
for all $\lambda \in (\delta^{\frac{1}{B}},1)$ and for all $K \geq 0.$ 
\end{lemma}
\begin{proof}
From PANDA's iterates we have that
\begin{align}
\mathbf{z}^{\perp}(k+1) = & \left(\boldsymbol{\Pi}^{\perp}_{\mathbf{1}_n}\mathbf{W}(k) \otimes \mathbf{I}_p\right)\mathbf{z}^{\perp}(k)  \\ & + \left(\boldsymbol{\Pi}_{\mathbf{1}_n}^{\perp}\otimes \mathbf{I}_p\right)\left(\mathbf{x}(k+1) - \mathbf{x}(k)\right), \nonumber
\end{align}
where we have used that $\boldsymbol{\Pi}_{\mathbf{1}_n}^{\perp}$ is idempotent and that $\mathbf{W}(k)$ is column stochastic. By applying the relation recursively we obtain
\begin{align}
\mathbf{z}^{\perp}(k+1) = & \left(\boldsymbol{\Pi}_{\mathbf{1}_n}^{\perp}\mathbf{W}_B(k) \otimes \mathbf{I}_p\right)\mathbf{z}^{\perp}(k+1-B) \\
&+\!\!\!\!\!\!\!\!\! \sum_{t=k-B+1}^k (\boldsymbol{\Pi}_{\mathbf{1}_n}^{\perp} \mathbf{W}_{k-t}(k) \otimes \mathbf{I}_p)(\mathbf{x}(t+1) - \mathbf{x}(t)). \nonumber
\end{align}
where we have used the fact that $\{\mathbf{W}(k)\}$ is column stochastic.
From here, by using Lemma \ref{lemma:mixing} and using that $\rho(\boldsymbol{\Pi}_{\mathbf{1}_n}^{\perp}\mathbf{W}(k)) \leq 1,$ we have
\begin{equation}
\|\boldsymbol{z}^{\perp}\| \leq \delta \|\mathbf{z}^{\perp}(k+1-B)\| + \sum_{t=1}^B\|\mathbf{x}(k+2-t) - \mathbf{x}(k+1-t)\|.
\end{equation}
Further, since $f^*$ has $\frac{1}{\mu}-$Lipschitz continuous gradients, which follows from Theorem \ref{theorem:dual}, we have that
\begin{align}
\label{eq:secondlemma}
\|\mathbf{z}^{\perp}(k+1)\| \leq & \delta \|\mathbf{z}^{\perp}(k+1-B)\|  \\
& + \sum_{t=1}^B \frac{1}{\mu}\|\mathbf{y}(k+1-t) - \mathbf{y}(k-t)\|.\nonumber
\end{align}
From here we can identify the terms in \eqref{eq:secondlemma} with the terms in \eqref{eq:afterlemma}. Hence, in order to establish this lemma we proceed as in the proof of Lemma \ref{lemma:a2} from equation \eqref{eq:afterlemma}.
\end{proof}
 In Lemma \ref{lemma:a5} we will establish that PANDA can be interpreted as a perturbed version of gradient descent applied on the dual problem. This interpretation allows us to use part of the framework in \cite{oracle} in order to establish our result. However, in order to increase the readability of the proof of Lemma \ref{lemma:a5} we will establish here an intermediate result.

\begin{lemma}[Gradient errors\label{lemma:gradient_errors}]
Given any $m>0$ strongly convex function $g$ with $K-$Lipschitz continuous gradients, a convex feasible set $\mathcal{Q},$ a sequence of vectors $\{\boldsymbol{\epsilon}(k)\}_{k=0}^{\infty},$ and the projected gradient descent method
\begin{align}
 \mathbf{x}(k+1):= \text{arg }\underset{\mathbf{x} \in \mathcal{Q}}{\text{min}} \quad \left( \nabla g(\mathbf{x}(k) + \boldsymbol{\epsilon}(k) \right)^T(\mathbf{x} - \mathbf{x}(k)) \nonumber\\
 \qquad \quad +  \frac{1}{2c}\|\mathbf{x} - \mathbf{x}(k)\|^2,
\end{align}
it holds that
\begin{align}
 \frac{m}{4}\|\mathbf{x} - \mathbf{x}(k)\|^2 \leq &  g(\mathbf{x}) - g_{\mathrm{err}}(\mathbf{x}(k),\boldsymbol{\epsilon}(k)) \\
& - \left( \nabla g(\mathbf{x}(k)) + \boldsymbol{\epsilon}(k) \right)^T(\mathbf{x} - \mathbf{x}(k)) \nonumber \\
\leq & K \|\mathbf{x}(k) - \mathbf{x}\|^2 + \frac{\|\boldsymbol{\epsilon}(k)\|^2}{2K},\,\forall \mathbf{x} \in \mathcal{Q}, \nonumber
\end{align}
where $g_{\mathrm{err}}(\mathbf{x},\boldsymbol{\epsilon}) \triangleq g(\mathbf{x}) - m\|\boldsymbol{\epsilon}(k)\|^2.$
\end{lemma}
\begin{proof}
The proof of this Lemma is very similar to the proof found in Section 2.3. of \cite{oracle} particularized to zero error in the evaluation of the function. The main difference being that the iterate number and error must be taken into account. The proof of this Lemma can be found for the case of time-varying errors in \cite{panda}. The proof is also included here for completeness.
Since the function $g$ as $K-$Lipschitz continuous gradients and is convex, we have that
\begin{align}
g(\mathbf{x}) & \leq g(\mathbf{x}(k)) - \boldsymbol{\epsilon}(k)^T(\mathbf{x} - \mathbf{x}(k)) \\
&+ \left( \nabla g(\mathbf{x}(k)) + \boldsymbol{\epsilon}(k) \right)^T(\mathbf{x} - \mathbf{x}(k)) + \frac{K}{2}\|\mathbf{x}(k) - \mathbf{x}\|^2 \leq \nonumber \\
& g(\mathbf{x}(k)) + \left( \nabla g(\mathbf{x}(k)) + \boldsymbol{\epsilon}(k) \right)^T(\mathbf{x}(k+1) - \mathbf{x}(k))  \nonumber\\
& + \left\| \boldsymbol{\epsilon}(k) \right\|\left\|\mathbf{x}(k) -\mathbf{x} \right\| + \frac{K}{2}\|\mathbf{x} - \mathbf{x}(k)\|^2,\, \forall \mathbf{x} \in \mathcal{Q}. \nonumber
\end{align}
Using the Peter-Paul inequality $\|\boldsymbol{\epsilon}(k)\|\|\mathbf{x}(k) - \mathbf{x}\| \leq \frac{\|\boldsymbol{\epsilon}(k)\|^2}{2K} + \frac{K}{2}\|\mathbf{x} - \mathbf{x}(k)\|^2$ yields
\begin{align}
\label{eq:lipschitz}
g(\mathbf{x}) - g(\mathbf{x}(k)) &\leq \\
& \left( \nabla g(\mathbf{x}(k)) + \boldsymbol{\epsilon}(k) \right)^T(\mathbf{x} - \mathbf{x}(k)) \nonumber \\
& + K\|\mathbf{x}(k) - \mathbf{x}\|^2 + \frac{\|\boldsymbol{\epsilon}(k)\|^2}{2K},\,\forall \mathbf{x} \in \mathcal{Q}. \nonumber
\end{align}
Analogously by using the fact that $g$ is $m-$strongly convex we have
\begin{align}
\label{eq:stronglyconvex}
g(\mathbf{x}) - g(\mathbf{x}(k)) & \geq \left( \nabla g(\mathbf{x}(k)) + \boldsymbol{\epsilon}(k) \right)^T(\mathbf{x} - \mathbf{x}(k)) \\
& \frac{m}{4}\|\mathbf{x} - \mathbf{x}(k)\|^2 - \frac{\|\boldsymbol{\epsilon}(k)\|^2}{m},\, \forall \mathbf{x} \in \mathcal{Q}. \nonumber
\end{align}
Combining both \eqref{eq:lipschitz} and \eqref{eq:stronglyconvex} yields
\begin{align}
\frac{m}{4}\|\mathbf{x} - \mathbf{x}(k)\|^2 & \leq g(\mathbf{x}) - g_{\mathrm{err}}(\mathbf{x}(k),\boldsymbol{\epsilon}(k)) \\
&- \left( \nabla g(\mathbf{x}(k)) + \boldsymbol{\epsilon}(k) \right)^T(\mathbf{x} - \mathbf{x}(k)) \nonumber \\
& \leq  K \| \mathbf{x}(k) - \mathbf{x} \|^2 + \frac{\|\boldsymbol{\epsilon}(k)\|^2}{2K},\, \forall \mathbf{x} \in \mathcal{Q},
\end{align}
which is the desired result.
\end{proof}
\begin{lemma}[(A5)\label{lemma:a5}] Under Assumptions \ref{assumption:mixing} and \ref{assumption:function}, $c \in \left(0,\frac{\mu}{2} \right],$ and $\lambda \in \left[\sqrt{1-\frac{c}{2L}},1 \right)$ it holds that
\begin{equation}
\|\mathbf{r}\|^{\lambda,K} \leq \sqrt{L\mu}\|\mathbf{z}\|^{\lambda,K} + 2 \|\mathbf{r}(0)\|.
\end{equation}
\end{lemma}
\begin{proof}
We start the analysis by equivalently re-writing PANDA as
\begin{subequations}
\begin{align}
 \mathbf{x}(k+1) := & \nabla f^*(\mathbf{y}(k)) \\
 \mathbf{z}(k+1) := & \mathbf{W}(k)\mathbf{z}(k) + \mathbf{x}(k+1) - \mathbf{x}(k) \\
 \mathbf{y}(k+1) := & \text{arg }\underset{\mathbf{y} \in \mathcal{Y}}{\text{min}} \left( \nabla f^*(\mathbf{y}(k)) + \boldsymbol{\epsilon}(k)\right)^T(\mathbf{y} - \mathbf{y}(k)) \nonumber \\
&+ \frac{1}{2c}\|\mathbf{y} - \mathbf{y}(k)\|^2, \label{eq:optimality_panda}
\end{align}
\end{subequations}
where $\mathcal{Y} \triangleq \{\mathbf{y} : (\boldsymbol{\Pi}_{\mathbf{1}_n}^{\perp} \otimes \mathbf{I}_p)\mathbf{y} = \mathbf{y}\}$ and $\boldsymbol{\epsilon}(k) = - \boldsymbol{\Pi}_{\mathbf{1}_n}^{\perp}\mathbf{z}(k+1).$
As can be seen by the equations above, PANDA can be interpreted as dual ascent in which the oracle providing the gradient information is inexact. Note that this can be done because the matrix $\boldsymbol{\Pi}_{\mathbf{1}_n}^{\perp}$ is idempotent. Further, since $f^*$ is strongly convex and has Lipschitz continuous gradients (cf. Theorem \ref{theorem:dual}) it follows from Lemma \ref{lemma:gradient_errors} that
\begin{align}
 \frac{1}{4L}\|\mathbf{y} - \mathbf{y}(k)\|^2 \leq &  f^*(\mathbf{y}) - f^*_{\text{err}}(\mathbf{y}(k),\boldsymbol{\epsilon}(k))^T(\mathbf{x} - \mathbf{x}(k)) \nonumber \\
 \leq &  \frac{1}{\mu}\|\mathbf{y}(k)-\mathbf{y}\|^2 + \frac{\mu\|\boldsymbol{\epsilon}(k)\|^2}{2}, \label{eq:combined}
\end{align}
where $f^{*}_{\text{err}}(\mathbf{y},\boldsymbol{\epsilon}) \triangleq f^*(\mathbf{y}) - \frac{1}{L}\|\boldsymbol{\epsilon}(k)\|^2.$ The proof that follows from here combines elements from the proof of Theorem 4 in \cite{oracle} and elements of the proof of Lemma 3.12 in \cite{interpret}.
For notational convenience let $r^2(k) \triangleq \|\mathbf{r}(k)\|^2.$ Then, for $k \geq 0$ it holds that
\begin{align}
r^2(k+1) = &\|\mathbf{y}(k+1) - \mathbf{y}^{\star} \|^2 = r^2(k) \nonumber\\
&+ 2(\mathbf{y}(k) - \mathbf{y}(k))^T(\mathbf{y}(k+1) - \mathbf{y}^{\star}) \nonumber \\
& - \|\mathbf{y}(k+1) - \mathbf{y}(k)\|^2. \label{eq:decrease}
\end{align}
By writing the optimality condition of \eqref{eq:optimality_panda} we obtain 
\begin{align}
\left( \nabla f^*(\mathbf{y}(k)) + \boldsymbol{\epsilon}(k) + \frac{1}{c}\left( \mathbf{y}(k+1) - \mathbf{y}(k) \right)\right)^T \\
\times \left(\mathbf{y} - \mathbf{y}(k+1)\right) \geq 0, \, \forall \mathbf{y} \in \mathcal{Y},
\end{align}
which can be particularized for $\mathbf{y}^{\star}$ yielding
\begin{align}
(\mathbf{y}(k+1) - \mathbf{y}(k))^T(\mathbf{y}(k+1) - \mathbf{y}^{\star}) \nonumber \\
\leq c(\nabla f^*(\mathbf{y}(k)) + \boldsymbol{\epsilon}(k))^T(\mathbf{y}^{\star} - \mathbf{y}(k+1)). \label{eq:particular}
\end{align}
By combining \eqref{eq:decrease} and \eqref{eq:particular} we obtain
\begin{align}
 r^2(k+1) \leq & r^2(k) + 2c(\nabla f^*(\mathbf{y}(k)) + \boldsymbol{\epsilon}(k))^T(\mathbf{y}^{\star} - \mathbf{y}(k)) \nonumber \\
&- 2c(\nabla f^*(\mathbf{y}(k)) - \boldsymbol{\epsilon}(k))^T(\mathbf{y}(k+1) - \mathbf{y}(k)) \nonumber \\
& + \frac{1}{2c}\|\mathbf{y}(k+1) - \mathbf{y}(k)\|^2. \label{eq:above}
\end{align}
By using \eqref{eq:combined} with $\mathbf{y} = \mathbf{y}(k+1)$ on \eqref{eq:above} we have
\begin{align}
r^2(k+1) \leq & r^2(k) + 2c \left( \nabla f^*(\mathbf{y}(k)) + \boldsymbol{\epsilon}(k) \right)^T(\mathbf{y}^{\star} - \mathbf{y}(k)) \nonumber \\
&\!\!\!\!\!\!\!\!\! -2c(f^*(\mathbf{y}(k+1)) - f^*_{\text{err}}(\mathbf{y}(k),\boldsymbol{\epsilon}(k))- \frac{\mu}{2}\|\boldsymbol{\epsilon}(k)\|^2),
\end{align}
as long as $c \in \left(0,\frac{\mu}{2} \right].$ By again using \eqref{eq:combined}, but this time with $\mathbf{y} = \mathbf{y}^{\star},$ we further obtain
\begin{align}
\label{eq:possibilities}
r^2(k+1) \leq & \left(1-\frac{c}{2L}\right)r^2(k) \\
& +2c (f^*(\mathbf{y}^{\star}) - f^*(\mathbf{y}(k+1))) + c \mu \|\boldsymbol{\epsilon}(k)\|^2. \nonumber
\end{align}
From here we do not proceed as in \cite{oracle} since we want to take advantage of the error's variability. We will instead proceed similarly to the final parts of the proof of Lemma 3.12 in \cite{interpret}. We consider two possibilities regarding $r^2(k+1).$ The first possibility being
\begin{equation}
\label{eq:possibilityA}
r^2(k+1) \geq L \mu \|\boldsymbol{\epsilon}(k)\|^2
\end{equation}
and the alternative
\begin{equation}
\label{eq:possibilityB}
r^2(k+1) < L \mu \|\boldsymbol{\epsilon}(k)\|^2.
\end{equation}
If \eqref{eq:possibilityA} occurs, we have
\begin{align}
\label{eq:possibility_A_consequence}
-2c (f^*(\mathbf{y}^{\star}) - f^*(\mathbf{y}(k+1))) \geq \frac{c}{L} \|\mathbf{y}(k+1) - \mathbf{y}^{\star}\|^2 = \nonumber \\
\frac{c}{L}r^2(k+1) \geq c \mu \|\boldsymbol{\epsilon}(k)\|,
\end{align}
which follows from the fact that $f^*$ is strongly convex and that $\nabla f^*(\mathbf{y}^{\star})^T(\mathbf{y} - \mathbf{y}^{\star}) = 0, \forall \mathbf{y} \in \mathcal{Y}.$ Note that this follows from the fact that $\nabla f^*(\mathbf{y}^{\star}) = \mathbf{x}^{\star}$ which lies in the orthogonal complement of $\mathcal{Y}.$ Combining \eqref{eq:possibility_A_consequence} with \eqref{eq:possibilities} yields
\begin{equation}
r^2(k+1) \leq \left(1 - \frac{c}{2L} \right)r^2(k).
\end{equation}
Hence, considering both possibilities \eqref{eq:possibilityA} and \eqref{eq:possibilityB} yields
\begin{equation}
r^2(k+1) \leq \text{max} \left\{ \left(1 - \frac{c}{2L} \right)r^2(k), L\mu \|\boldsymbol{\epsilon}(k)\|^2\right\}. \label{eq:recursive}
\end{equation}
Applying \eqref{eq:recursive} recursively we obtain
\begin{align}
r^2(k+1) \leq & \text{max} \left\{ \left(1 - \frac{c}{2L} \right)^{k+1}r^2(0), \right. \\
& \left. L \mu \underset{t=0,\hdots,k}{\text{sup}} \left\{ \left(1 - \frac{c}{2L} \right)^t \|\boldsymbol{\epsilon}(k-t)\|^2\right\} \right\}. \nonumber
\end{align}
By taking square roots on both sides we obtain
\begin{align}
r(k+1) \leq &\left(\sqrt{1 - \frac{c}{2L}}\right)^{k+1}r(0) \\
&+\sqrt{L \mu}\underset{t=0,\hdots,k}{\text{sup}} \left\{ \left(\sqrt{1 - \frac{c}{2L} }\right)^t\| \boldsymbol{\epsilon}(k-t) \| \right\}. \nonumber
\end{align}
Now select $\lambda$ such that $\gamma \triangleq (\lambda)^{-2}(1-\frac{c}{2L}) \leq 1.$ Then, we have
\begin{align}
\lambda^{-(k+1)}r(k+1) \leq & (\sqrt{\gamma})^{k+1}r(0) \\
&\!\!\!\!\!\!\!\!\!\!\! + \lambda^{-1}\sqrt{L \mu} \underset{t=0,\hdots,k}{\text{sup}} \left\{ \lambda^{-(k-t)}\sqrt{\gamma ^t}\|\boldsymbol{\epsilon}(k-t)\| \right\} \nonumber\\
& \!\!\!\!\!\!\!\!\!\!\leq r(0) + \lambda^{-1}\sqrt{L \mu}\underset{t=0,\hdots,k}{\text{sup}} \left\{ \lambda^{-t}\|\boldsymbol{\epsilon}(t)\| \right\}.\nonumber
\end{align}
By taking $\underset{k=0,\hdots,K}{\text{sup}}$ on both sides we obtain
\begin{equation}
|r|^{\lambda,K} \leq 2r(0) + \lambda^{-1}\sqrt{L \mu}\|\boldsymbol{\epsilon}\|^{\lambda,K}.
\end{equation}
In order to establish that the last arrow in \eqref{eq:cycle} it suffices to recall that $\|\boldsymbol{\epsilon}(k)\| = \|\mathbf{z}^{\perp}(k+1)\|$ and that $r(k) = \|\mathbf{r}(k)\|.$ Therefore, we conclude
\begin{equation}
\|r\|^{\lambda,K} \leq 2r(0) + \sqrt{L \mu}\|\mathbf{z}^{\perp}\|^{\lambda,K},
\end{equation}
for $\lambda \in [\sqrt{1-\frac{c}{2L}},1)$ and $c \in (0,\frac{\mu}{2}].$
\end{proof}
\subsection{Applying the small gain theorem \label{subsection:small}}
We are now ready to prove Theorem \ref{theorem:panda}. For this we use the small gain theorem (c.f. Theorem \ref{theorem:gain}). In order to establish convergence via the small gain theorem we have to find a $\lambda$ and a $c$ such that
\begin{subequations}
\label{eq:conditions_small}
\begin{align}
\frac{(1-\lambda^B)^2c}{(1-\lambda)^2(\lambda^B - \delta)^2} < \frac{\mu}{2\sqrt{\kappa}
}, \label{eq:conditions_small1}\\
0 < c \leq \frac{\mu}{2}, \label{eq:conditions_small2}\\
\sqrt{1-\frac{c}{2L}} \leq \lambda \label{eq:conditions_small3}\\
\lambda \in \left( \delta^{1/B},1 \right). \label{eq:conditions_small4}
\end{align}
\end{subequations}
where $\kappa \triangleq \frac{L}{\mu}.$
Condition \eqref{eq:conditions_small1} follows from requiring that $\gamma_1\hdots\gamma_5 < 1$ (c.f. Section \ref{subsection:proofstructure}) hold. Conditions \eqref{eq:conditions_small2},\eqref{eq:conditions_small3} are required to guarantee that the bound established in Lemma \ref{lemma:a5} holds. Finally, \eqref{eq:conditions_small4} is required so that the bounds established in Lemmas \ref{lemma:a2} and \ref{lemma:a4} hold.
Before proceeding with the proof we require the following lemma.
\begin{lemma}\label{lemma:bound}
For $\lambda \in [\frac{\sqrt{3}}{2}, 1)$ and any finite $B$ it holds that
\begin{equation}
\label{eq:lemma_bound}
\frac{(1-\lambda^B)^2}{(1-\lambda)^2} \leq 1.
\end{equation}
\end{lemma}
\begin{proof}
In order to first establish the bound we establish that within the range $\lambda \in [0,1)$ the function is decreasing in $\lambda.$ This can be done by seeing that the inequality holds with equality if $B=1.$ However, when $B>1$ the numerator is always smaller than the denominator leading implying strict inequality.
\end{proof}
Note that since $c \leq \frac{\mu}{2}$ and $\lambda \geq \sqrt{1-\frac{c}{2L}}$ we have that
$\lambda \geq \sqrt{1 - \frac{\mu}{4L}} \geq \frac{\sqrt{3}}{2}$ and hence Lemma \ref{lemma:bound} holds for all valid values of $\lambda.$
By using Lemma \ref{lemma:bound} any pair $\lambda,\,c$ fulfilling
\begin{subequations}
\label{eq:stronger_conditions}
\begin{align}
c \leq \frac{\mu(\lambda^B - \delta)^2}{2\sqrt{\kappa}} \label{eq:c1}\\
c \geq 2L(1-\lambda^2) \\
\lambda \in \left(\delta^{1/B},1 \right) \\
c \leq \frac{\mu}{2} \label{eq:c2}
\end{align}
\end{subequations}
fulfils \eqref{eq:conditions_small}. Any $c$  fulfilling \eqref{eq:c1} will fulfil \eqref{eq:c2}. In other words, we are now to find a pair $(\lambda,c)$ that fulfil
\begin{equation}
c \in \left[2L(1-\lambda^2), \frac{\mu(\lambda^B - \delta)^2}{2\sqrt{\kappa}}\right],
\end{equation}
with $\lambda \in (\delta^{1/B},1).$ In order to be able to find such a pair in closed form, we restrict the lower bound on $c$ further by instead requiring
\begin{equation}
\label{eq:cbound}
c \in \left[2L(1-\lambda^{2B}), \frac{\mu(\lambda^B - \delta)^2}{2\sqrt{\kappa}}\right].
\end{equation}
Note that the lower bound in \eqref{eq:cbound} is a monotonically decreasing function of $\lambda$ while the upper bound is a monotonically increasing function of $\lambda.$ We will now establish that the upper and lower bounds intersect at a point $\lambda_{\bar{c}}$ where $\lambda_{\bar{c}} \in (\delta^{1/B},1).$ For this we find $\lambda_{\bar{c}}^B$ that solves
\begin{equation}
\left(\frac{1}{2}+2\kappa^{\frac{3}{2}}\right) \lambda^{2B} - \delta \lambda^B + \frac{\delta^2}{2} -2\kappa^{\frac{3}{2}} = 0
\end{equation}
which is given by
\begin{equation}
\label{eq:lambda_sol}
\lambda_{\bar{c}}^B = \frac{ \delta + \sqrt{\delta^2 + (4\kappa^{\frac{3}{2}}-\delta^2)(1+4\kappa^{\frac{3}{2}})}}{1+4\kappa^{\frac{3}{2}}}.
\end{equation}
$\lambda_{\bar{c}}^B$ can be easily verified to be both strictly larger than $\delta$ and strictly smaller than 1. In order to verify that it is smaller than 1 it is sufficient to show that the derivative is strictly positive for $\delta \in \left[0,1 \right]$ implying that the function \eqref{eq:lambda_sol} is monotonically increasing in $\delta.$ By evaluating \eqref{eq:lambda_sol} at $\delta = 1$ we see that $\lambda_{\bar{c}}^B < 1$ for $\delta \geq 1.$ Further, it can be seen that $\lambda_{\bar{c}}^B > \delta$ by algebraic manipulation of the inequality $\lambda{\bar{c}}^B > \delta$ and seeing that the resulting expression will hold true as long as $\delta < 1.$

 When both the upper and lower bound in \eqref{eq:cbound} intersect, the viable values of $c$ collapse to a single point $\bar{c}.$ This point can be found by evaluating $\lambda_{\bar{c}}^B$ in either the lower or higher bound of \eqref{eq:cbound}. For ease of computation we choose the upper bound 
\begin{align}
&\bar{c} = \frac{\mu}{2}\left(\frac{16\kappa^{\frac{5}{2}}-4\kappa(1-\delta^2)}{(1+4\kappa^{\frac{3}{2}})^2} \right).
\end{align}
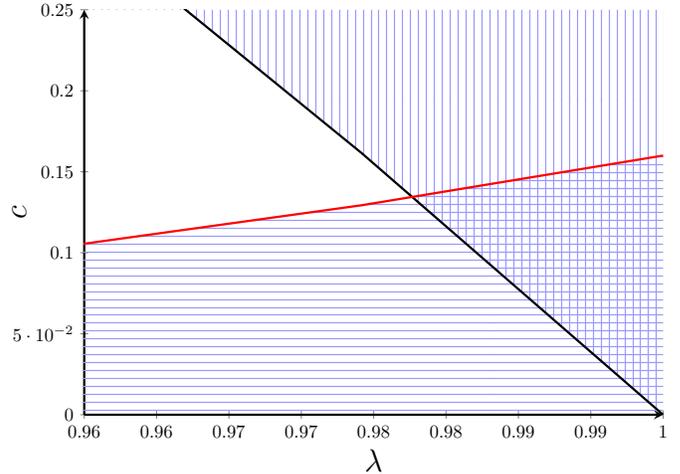
\begin{figure}
\scalebox{0.7}{
\begin{tikzpicture}
\begin{axis}[
width=4.951in,
height=3.653in,
at={(0.831in,0.513in)},
axis line style=very thick,
axis x line=bottom,
axis y line=left,
	ymin=0, ymax=0.25, xmin=0.96, xmax=1,
	xlabel style = {font = \color{black} \LARGE},
	ylabel style = {font = \color{black} \LARGE},
	xlabel=$\lambda$, ylabel=$c$,
]
\addplot [draw=none,
pattern=vertical lines, pattern color=blue!40,
 domain=0.5:1]{2 * (1 - x^4)} \drawge ;
\addplot [draw=none,
 pattern=horizontal lines, pattern color=blue!40,
 domain=0.5:1]{(x^4-0.2)^2/4} \drawle;
 \addplot[very thick, domain=0.5:1]{2 * (1 - x^4)};
 \addplot[very thick, red, domain=0.5:1]{(x^4-0.2)^2/4};
\end{axis}
\end{tikzpicture}}
\caption{Upper (red) and lower bounds (black) of $c$ as a function of $\lambda$ (c.f. \eqref{eq:cbound}) with $B = 2, \delta = 0.2,$ $L = \mu = \kappa =1.$ $(\lambda_{\bar{c}},\bar{c})$ is the intersection between the two curves. Viable options for $(\lambda_,c)$ are those in the region shaded both horizontally and vertically.\label{fig:intersection}}
\end{figure}

Recall now that the upper and lower bounds in \eqref{eq:cbound} are monotonically increasing and decreasing functions of $\lambda$ respectively. We have already established that the two curves intersect at a point $\lambda_{\bar{c}}$ with $\lambda_{\alpha} > \delta^{1/B}.$ Therefore, the feasible values of $c$ are found in between the two curves after they have intersected. This is clarified further in Fig. \ref{fig:intersection}.
Note that since the lower bound in \eqref{eq:cbound} is monotonically decreasing in $\lambda$ as the values of $c$ become smaller the smallest possible $\lambda$ to chose is that which will match the lower bound. Hence,
\begin{equation}
\text{if } c\in \left(0,\bar{c} \right], \, \lambda = \sqrt[2B]{1-\frac{c}{2L}}. 
\end{equation}  
On the other hand, whenever $c$ becomes large the smallest possible $\lambda$ to chose is that matching the upper bound, i.e.
\begin{equation}
\text{if } c \in \left(\bar{c},\frac{\mu (1-\delta)^2}{2\sqrt{\kappa}}\right), \, \lambda = \sqrt[B]{\delta + \sqrt{\frac{2c\sqrt{\kappa}}{\mu}}}
\end{equation}
where the upper bound on $c$ comes from setting $\lambda = 1.$
Hence, since we have established the arrows and the existence of an interval of $(\lambda,c)$ that fulfills the requirements of the small gain theorem, we conclude that Theorem 1 has been proven.
%\section{Directed Networks}
%\subsection{Algorithm and Intuition}
%\begin{algorithm}
%\caption{Push-PANDA \label{alg:push_panda}}
%\begin{algorithmic}[1]
%\State Choose step size $c > 0$ and pick $\mathbf{z}(0)= \mathbf{x}(0) = \mathbf{0},$ $\mathbf{v}(0) = \mathbf{1}_n$ and $\mathbf{y}(0)$ such that $(\boldsymbol{\Pi}_{\mathbf{1}_n} \otimes \mathbf{I}_p)\mathbf{y}(0) = \mathbf{0}.$
%\For{$k=0,1,\hdots$} Each agent $i$:
%\State Computes
%\begin{equation}
%\mathbf{x}_i(k+1) :=  \underset{\mathbf{x}_i \in \mathbb{R}^{p}}{\text{min}}  f_i(\mathbf{x}_i) - \mathbf{y}_i(k)^T\mathbf{x}_i
%\end{equation}
%\State Exchanges $v_i(k)$ with $j \in \mathcal{N}_i^{\text{out}}(k)$
%\State Computes
%\begin{equation}
%v_i(k+1) : = C_{ii}(k)v_i(k) + \sum_{j \in \mathcal{N}_i^{\text{in}}(k)}C_{ij}(k)v_j(k)
%\end{equation}
%\State Exchanges $\mathbf{z}_i(k)$
% with $\mathcal{N}_i^{\text{out}}(k).$ 
%\State Computes 
%\begin{align}
%\mathbf{z}_i(k+1) :=  C_{ii}(k)\mathbf{z}_i(k) + \sum_{j \in \mathcal{N}_i^{\text{in}}(k)}C_{ij}(k)\mathbf{z}_j(k) \\+ \left(\mathbf{x}_i(k+1) - \mathbf{x}_i(k)\right) \nonumber
%\end{align}
%\State Exchanges $\mathbf{x}_i(k+1) - \frac{1}{v_i(k+1)}\mathbf{z}_i(k+1)$ with $\mathcal{N}_i^{\text{out}}(k)$
%\State Computes 
%\begin{align}\mathbf{y}_i(k+1) := \mathbf{y}_i(k) - c\left(\mathbf{x}_i(k+1) -\frac{1}{v_i(k+1)}\mathbf{z}_i(k+1)\right) \nonumber \\
%+ c \sum_{j \in \mathcal{N}_i^{\text{in}}} C_{ij}(k)\left(\mathbf{x}_j(k+1) - \frac{1}{v_i(k+1)}\mathbf{z}_j(k+1)\right) 
%\end{align}
%\EndFor
%\end{algorithmic}
%\end{algorithm}
%\subsection{R-linear convergence}
\section{Numerical Experiments \label{section:numerical}}
Consider the following decentralized estimation problem. Each agent $i$ obtains the measurement $\mathbf{b}_i \in \mathbb{R}^{d}$ via a measurement matrix $\mathbf{H}_i \in \mathbb{R}^{d \times p}.$ The goal of the set of agents is to collaboratively solve the ridge regression problem \cite{time_nesterov}
\begin{equation}
\label{eq:numerical_centralized}
\underset{\mathbf{x} \in \mathbb{R}^{p}}{\min} \quad \frac{1}{2nd}\sum_{i=1}^n\left\|\mathbf{H}_i\mathbf{x} - \mathbf{b}_i\right\|^2 + \frac{r}{2}\|\mathbf{x}\|^2
\end{equation}
over a time-varying network represented by the graph $\mathcal{G}(\mathcal{V},\mathcal{E}(k)).$
The elements of the vector  $\mathbf{x}^{\star}$ which solves \eqref{eq:numerical_centralized} are generated independently and according to the distribution  $ \mathcal{N}(0,10).$ Then, $\mathbf{b}_i = \mathbf{H}_i \mathbf{x} + \boldsymbol{\epsilon}_i,$ where the elements of $\boldsymbol{\epsilon}_i$ are generated independently and according to the distribution $\mathcal{N}(0,0.1).$ Similarly, the elements of the matrices $\mathbf{H}_i \in \mathbb{R}^{d \times p}$ are generated independently and according to the distribution $\mathcal{N}(0,0.1).$ Further, the sequence of graphs $\mathcal{G}(\mathcal{V},\mathcal{E}(k))$ is randomly generated randomly. In particular, each bi-directional arc will belong to the set of edges with probability $\pi,$ i.e. $P((i,j) \in \mathcal{E}(k)) = \pi,$ and $P((i,j) \in \mathcal{E}(k)| (j,i) \in \mathcal{E}(k)) = 1.$ Also, $P((i,j) \in \mathcal{E}(k),(k,l) \in \mathcal{E}(k)) = P((i,j) \in \mathcal{E}(k))P((k,l) \in \mathcal{E}(k)),$  where $(k,l) \neq (i,j)$ and $(k,l) \neq(j,i).$ All experiments are performed over a network of 10 nodes, i.e. $n = 10,$ with $d = 3$ and $p = 5.$ The Metropolis Hastings matrix is chosen for both PANDA and DIGing, i.e. the mixing matrix for both methods is given by 
\begin{equation}
[W(k)]_{i,j} = \begin{cases}
\frac{1}{\max\{d_i(k),d_j(k)\}} & \text{if } (i,j) \in \mathcal{E}(k) \\
1 - \sum_{j \in \mathcal{N}_{i}(k)} w_{ij}(k) & \text{if } i = j \\
0 & \text{otherwise}
\end{cases},
\end{equation}
while the matrix used for Dual Ascent is the graph's Laplacian matrix at each iteration.
\begin{figure}
\scalebox{0.57}{\input{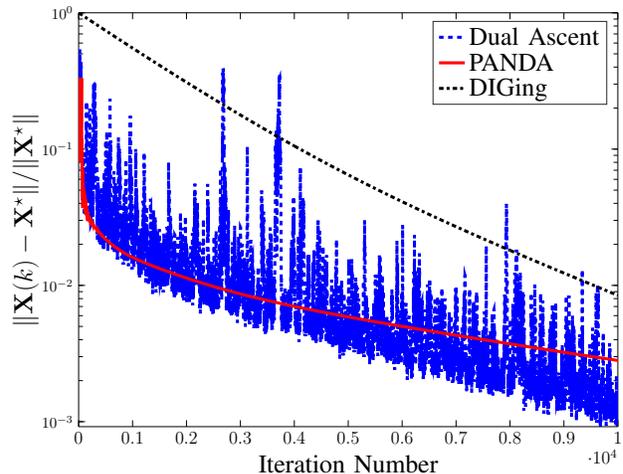}}
\caption{\label{fig:low} Ridge regression problem with 10 nodes over sparsely connected networks, $p = 3, $ $d=5$ and $r = 0.001.$ The network is generated using $\pi = 0.1.$ The condition number $\kappa = 4.3 \cdot 10^3,$ is computed a posteriori and the step-sizes are $\alpha_{dd}=9 \times 10^{-5},$ $c = 5 \times 10^{-5},$ $\alpha_{ing} = 0.06,$ for dual decomposition, PANDA and DIGing respectively. }
\end{figure}
\begin{figure}
\scalebox{0.57}{\input{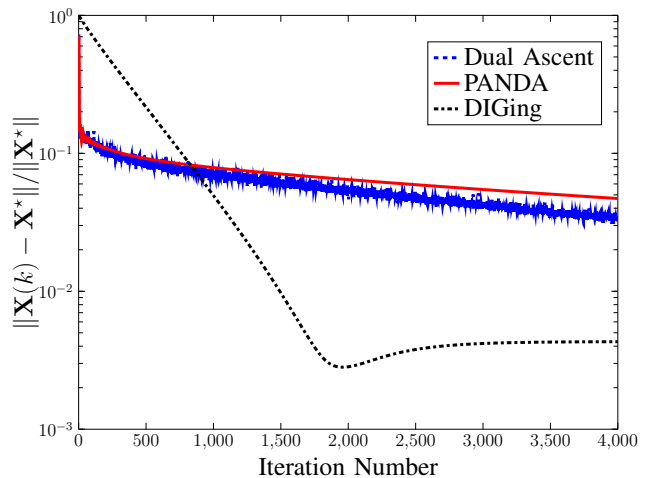}}
\caption{\label{fig:bad} Badly conditioned ridge regression problem with 10 nodes, $p = 3$, $d=5,$ and $r = 10^{-4}.$ The network is generated using $\pi=0.5.$ The condition number  $\kappa = 4.09 \times 10^{4}$ is computed a posteriori and the step-sizes are $\alpha_{dd} = 3 \times 10^{-6},$ $c = 9 \times 10^{-6},$ $\alpha_{ing} = 0.09,$ for dual decomposition, PANDA and DIGing respectively.}
\end{figure}
\begin{figure}
\scalebox{0.57}{\input{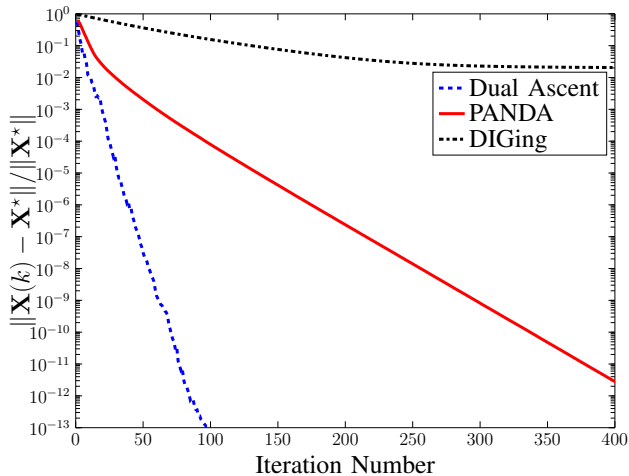}}
\caption{\label{fig:good} Well conditioned ridge regression problem with 10 nodes, $p = 3,$ $d=5,$ and $r = 1.$ The network is generated using $\pi = 0.5.$ The condition number $\kappa = 3.9$ is computed a posteriori and the step-sizes are $\alpha_{dd} = 0.03 ,$ $c = 0.03 ,$ $\alpha_{ing} = 0.09,$ for dual decomposition, PANDA and DIGing respectively.}
\end{figure}
 Figures \ref{fig:low}-\ref{fig:good} illustrate the performance of DIGing \cite{interpret}, PANDA and Dual Decomposition \cite{dual_decomposition,fenchel}. The performance is measured in terms of the quantity $\frac{\|\mathbf{X}(k) -\mathbf{X}^{\star}\|_{\text{F}}}{\|\mathbf{X}^{\star}\|_{\text{F}}},$
  where $\mathbf{X}(k) \triangleq [\mathbf{x}_1(k)^{\text{T}}; \hdots ; \mathbf{x}_N(k)^{\text{T}}]$ and
  $\mathbf{X}^{\star} \triangleq [\mathbf{x}^{\star \text{T}}; \hdots; \mathbf{x}^{\star \text{T}}].$ The step-size for each of the methods has been optimized by hand to yield the best convergence rate. 
 
 From the numerical analysis we can conclude that which method performs best depends on the instance of the problem \eqref{eq:numerical_centralized}. However, as Fig. \ref{fig:low} indicates, when the network connectivity becomes worse, i.e., the probability of the nodes being connected decreases, PANDA outperforms both Dual Decomposition and DIGing. Recall that PANDA comes at no added computational expense when compared to Dual Decomposition, and is cheaper in terms of information exchange than DIGing.
Perhaps unsurprisingly, PANDA has some of the advantages and disadvantages of both Dual Decomposition and DIGing. In particular, when the objective function is well conditioned, i.e. $\kappa \not \gg 1,$ both Dual Decomposition and PANDA will perform better than DIGing as seen in Fig. \ref{fig:good} while PANDA performs worse than Dual Decomposition. On the other hand of the spectrum, when the condition number becomes very large, as seen in Fig. \ref{fig:bad} Dual Decomposition and PANDA exhibit practically identical performance while the iterates of PANDA appear to be more stable. In this case, DIGing performs better than both Dual Decomposition and PANDA.
\section{Conclusion}
In this paper we proposed PANDA, a dual ascent based
method for time varying graphs. We established that it convergences
R-linearly over time varying graphs and hence is the first dual method to provably do so. 
Finally, we compare its numerical performance to that of Dual Decomposition and DIGing. One of the advantages both PANDA and Dual Decomposition to DIGing is that they
require communicating half as many quantities as DIGing per
iteration. On the other hand, both PANDA’s and Dual Decomposition's iterates are in general computationally more expensive than those of DIGing. The extend to which the computations are more expensive will exclusively depend on the problem. This
implies that while DIGing may be more suitable for scenarios
in which the computation cost is high PANDA may be more
suitable for scenarios in which the communication cost is
high. Finally, we see that for sparsely connected networks PANDA outperforms Dual Decomposition. 
\appendices
%\section*{Acknowledgment}

\ifCLASSOPTIONcaptionsoff
  \newpage
\fi
%\begin{thebibliography}{1}
\bibliographystyle{IEEEbib.bst}
\bibliography{bibliography}
%\end{thebibliography}
\begin{IEEEbiography}
[{\includegraphics[width=1in,height=1.25in,clip,keepaspectratio]{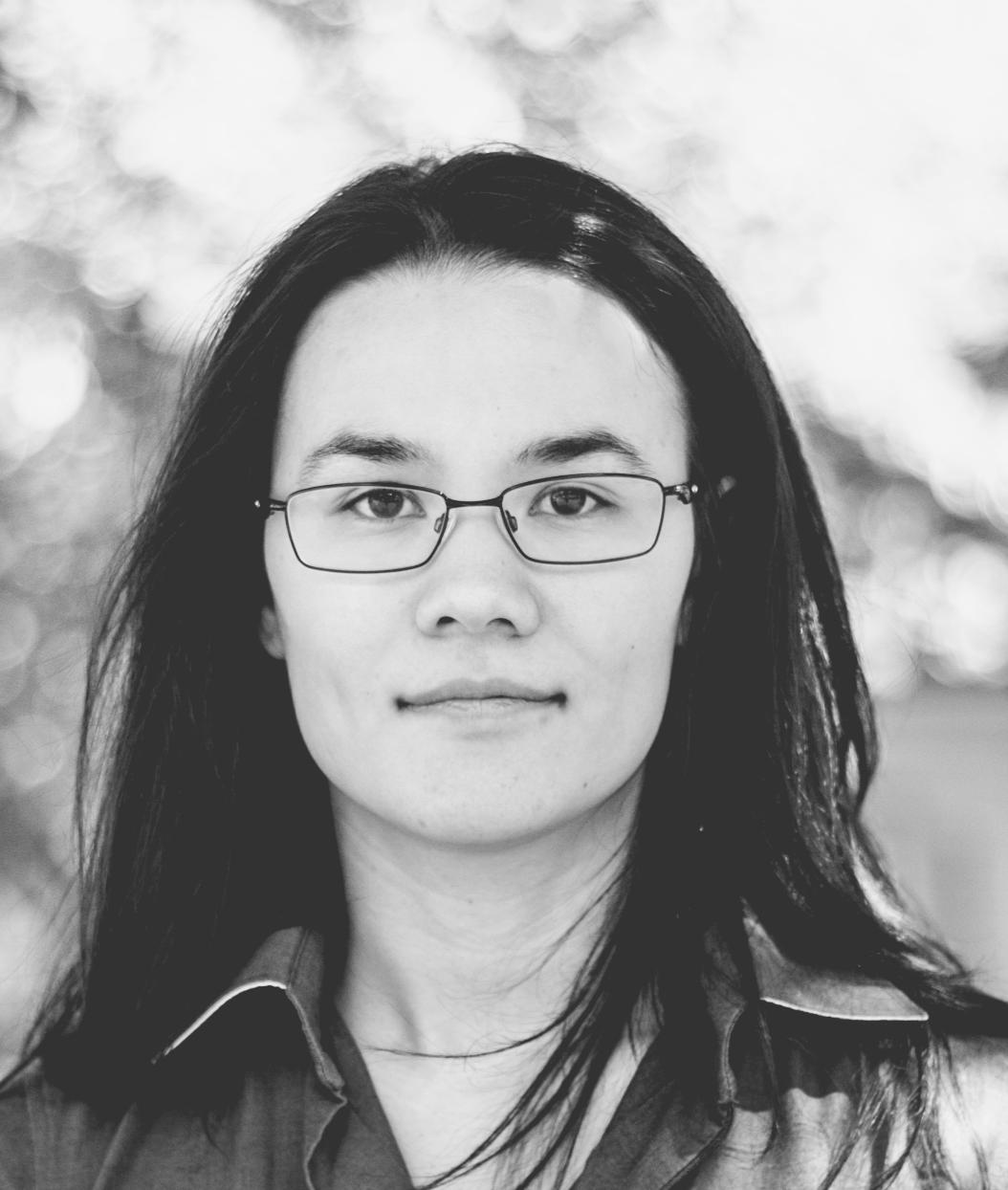}}]{Marie Maros} (S'14) received the M.Sc. in electrical engineering from the Technical University of Catalonia (UPC), Barcelona, Spain, and from the Royal Institute of Technology (KTH), Stockholm, Sweden in 2014. In August 2014 she joined the Information Science and Engineering Lab within the School of Electrical Engineering at KTH, Stockholm, Sweden where she is currently pursuing her Ph.D. degree in electrical engineering. Her research interests include distributed optimization for time varying systems.
\end{IEEEbiography}
\begin{IEEEbiography}[{\includegraphics[width=1in,height=1.25in,clip,keepaspectratio]{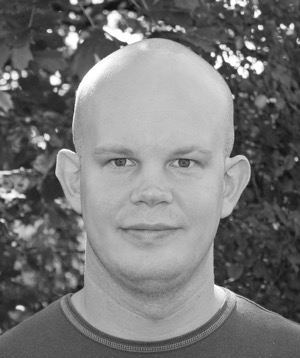}}]{Joakim Jal\'{e}n} (S'03-M'08-S'13) received the M.Sc.\ and Ph.D.\ in electrical engineering from the Royal Institute of Technology (KTH), Stockholm, Sweden in 2002 and 2007 respectively. From July 2007 to June 2009 he held a post-doctoral research position at the Vienna University of Technology, Vienna, Austria. He also studied at Stanford University, CA, USA, from September 2000 to May 2002, and worked at ETH, Z\"{u}rich, Switzerland, as a visiting researcher, from August to September, 2008. In July 2009 he joined the Information Science and Engineering Lab within the School of Electrical Engineering at KTH, Stockholm, Sweden, as an Assistant Professor. He was an associate editor for the IEEE Communications Letters between 2009 and 2011, an associate editor for the IEEE Transactions in Signal Processing between 2012 and 2016, and is a member of the IEEE Signal Processing for Communications and Networking Technical Committee (SPCOM-TC). 

For his work on MIMO communications, Joakim has been awarded the IEEE Signal Processing Society’s 2006 Young Author Best Paper Award, the Distinguished Achievement Award of NEWCOM++ Network of Excellence in Telecommunications 2007–2011, and the best student conference paper award at IEEE ICASSP 2007. He is also a recipient of the Ingvar Carlsson Career Award issued in 2009 by the Swedish Foundation for Strategic Research. His recent work includes work on signal processing for biomedical data analysis, and the automated tracking of (biological) cell migration and morphology in time-lapse microscopy in particular. Early work in this field was awarded a conference best paper award at IEEE ISBI 2012, and subsequent work by the group has been awarded several Bitplane Awards in connection to the ISBI cell tracking challenges between 2013 and 2015.
\end{IEEEbiography}
%\vfill
%\newpage
%\vfill
\end{document}